\documentclass[12pt]{amsart}
\usepackage{amsfonts}
\usepackage{umlaut,german,graphicx}
\usepackage{amssymb}
\newcommand{\N}{\mathbb{N}}
\newcommand{\Z}{\mathbb{Z}}
\newcommand{\C}{\mathbb{C}}

\newcommand{\g}{\mathfrak{g}}
\newcommand{\h}{\mathfrak{h}}

\newcommand{\p}{\mathfrak{p}}

\usepackage{amsthm}
\usepackage{amsmath}
\usepackage{amscd}

\usepackage{amstext}

\newtheorem{theorem}{Theorem}
\newtheorem{lemma}{Lemma}
\newtheorem{proposition}{Proposition}
\newtheorem{corollary}{Corollary}
\newtheorem{definition}{Definition}

\newcommand{\on}[2]{\setbox0=\hbox{$#1$}\setbox1=\hbox{$#2$}%
            \dimen0=\wd0\advance\dimen0 by \wd1\divide\dimen0 by 2%
             \ifdim\wd0>\wd1$#1$\hskip-\dimen0$#2$\advance\dimen0 by -\wd1%
              \else$#2$\hskip-\dimen0$#1$\advance\dimen0 by -\wd0%
             \fi%
            \hskip\dimen0}
\newcommand{\onn}[2]{\makebox{\on{#1}{#2}}}
\newcommand{\cross}{\smash\times\vphantom\bullet}
\newcommand{\leftblob}[1]{\hspace*{-1ex}\onn{\raisebox{-.5ex}
     {$\genfrac{}{}{0pt}{}{\hfill\hrulefill}{\hphantom{\hspace*{2ex}#1}}$}}
                {\stackrel{#1}{\bullet}}}
\newcommand{\rightblob}[1]{\onn{\raisebox{-.5ex}
     {$\genfrac{}{}{0pt}{}{\hrulefill\hfill}{\hphantom{\hspace*{2ex}#1}}$}}
                 {\stackrel{#1}{\bullet}}\hspace*{-1ex}}
\newcommand{\middleblob}[1]{\onn{\raisebox{-.5ex}
     {$\genfrac{}{}{0pt}{}{\hrulefill}{\hphantom{\hspace*{2ex}#1}}$}}
                 {\stackrel{#1}{\bullet}}}
\newcommand{\leftcross}[1]{\hspace*{-1ex}\onn{\raisebox{-.5ex}
     {$\genfrac{}{}{0pt}{}{\hfill\hrulefill}{\hphantom{\hspace*{2ex}#1}}$}}
                {\stackrel{#1}{\cross}}}

\newcommand{\oo}[2]{\leftblob{#1}\hspace*{-1ex}\rightblob{#2}}

\newcommand{\xo}[2]{\leftcross{#1}\hspace*{-1ex}
               \rightblob{#2}}       
\newcommand{\ooo}[3]{\leftblob{#1}\hspace*{-1ex}
               \middleblob{#2}\hspace*{-1ex}\rightblob{#3}}
\newcommand{\xoo}[3]{\leftcross{#1}\hspace*{-1ex}
               \middleblob{#2}\hspace*{-1ex}\rightblob{#3}}

\newcommand{\oooo}[4]{\leftblob{#1}\hspace*{-1ex}\middleblob{#2}\hspace*{-1ex}
               \middleblob{#3}\hspace*{-1 ex}\rightblob{#4}}

\newcommand{\semioverbrace}[1]{\mathop{\vbox{\ialign{##\crcr\noalign{\kern3pt}
      \semidownbracefill\crcr\noalign{\kern3pt\nointerlineskip}
      $\hfil\displaystyle{#1}\hfil$\crcr}}}\limits}
\newcommand{\arrowhead}{
                \begin{picture}(0,0)\put(3,.605){\vector(1,0){0}}\end{picture}}
              
\newcommand{\xooo}[4]{\leftcross{#1}\hspace*{-1ex}
               \middleblob{#2}\hspace*{-1 ex}\middleblob{#3}\hspace*{-1ex}\rightblob{#4}}

\makeatletter
\newcommand{\semidownbracefill}{$\m@th\braceld\leaders\vrule\hfill\braceru
  \bracelu\leaders\vrule\hfill\arrowhead$}


\begin{document}
\title{Invariant differential pairings}
\author{Jens Kroeske}
\address{School of Mathematical Sciences, University of Adelaide,
SA 5005, Australia}
\email{jens.kroeske@student.adelaide.edu.au}
\renewcommand{\subjclassname}{\textup{2000} Mathematics Subject Classification}
\maketitle

\begin{abstract}
In this paper the notion of an $M$-th order invariant bilinear differential pairing is introduced and a formal definition is given. If the manifold has an AHS structure, then various first order pairings are constructed. This yields a classification of all first order invariant bilinear differential pairings on homogeneous spaces with an AHS structure except for certain totally degenerate cases. Moreover higher order invariant bilinear differential pairings are constructed on these homogeneous spaces which leads to a classification on complex projective space for the non degenerate cases. A degenerate case corresponds to the existence of an invariant linear differential operator.
\end{abstract}

\section{Introduction}

It is generally known (see ~\cite{pr},  p.~202), that on an arbitrary manifold $\mathcal{M}$ one can write down the Lie derivative $\mathcal{L}_{X}\omega_{b}$ of a one-form $\omega_{b}\in\Omega^{1}(\mathcal{M})$ with respect to a vector field $X^{a}\in T\mathcal{M}$ in terms of an arbitrary torsion-free connection $\nabla_{a}$ as
$$\mathcal{L}_{X}\omega_{b}=X^{a}\nabla_{a}\omega_{b}+\omega_{a}\nabla_{b}X^{a},$$
where the indices used are abstract in the sense of ~~\cite{pr}. This {\bf pairing} is obviously linear in $X^{a}$ and in $\omega_{b}$, i.e.~{\bf bilinear}, {\bf first order} and
{\bf invariant} in the sense that it does not depend upon a specific choice of connection. One can specify an equivalence class of connections and ask for invariance under change of connection within this equivalence class. In conformal geometry, for example, one deals with an equivalence class of connections that consists of the Levi-Civita connections that correspond to metrics in the conformal class. This paper will mainly deal with {\bf projective} {\bf geometry} and the projective equivalence class of connections consists of those torsion-free connections which induce the same (unparameterised) geodesics. This is equivalent (see ~\cite{e}, proposition~1) to saying that $\nabla_{a}$ and $\hat{\nabla}_{a}$ are in the same equivalence class
if and only if there is a one form $\Upsilon_{a}$, such that
$$\hat{\nabla}_{a}\omega_{b}=\nabla_{a}\omega_{b}-2\Upsilon_{(a}\omega_{b)},$$
where round brackets around indices denoted symmetrisation, i.e. $\Upsilon_{(a}\omega_{b)}=\frac{1}{2}(\Upsilon_{a}\omega_{b}+\Upsilon_{b}\omega_{a})$. As a consequence, the difference between the two connections when acting on sections of any weighted tensor bundle can be deduced (see ~\cite{e}) and the invariance of any expression can be checked by hand. For vector fields, for example, we have
$$\hat{\nabla}_{b}X^{a}=\nabla_{b}X^{a}+\Upsilon_{b}X^{a}+\Upsilon_{c}X^{c}\delta_{b}{}^{a},$$
so the invariance of the Lie derivative can be checked directly.
It is also clear that
$$X^{a}\nabla_{[a}\omega_{b]}$$
is a first order bilinear invariant differential pairing, where square brackets around indices denote skewing, i.e. $\nabla_{[a}\omega_{b]}=\frac{1}{2}(\nabla_{a}\omega_{b}-\nabla_{b}\omega_{a}$). Therefore there are at least two first order bilinear invariant pairings
$$T\mathcal{M}\times\Omega^{1}(\mathcal{M})\rightarrow\Omega^{1}(\mathcal{M}).$$
\par
To obtain a more interesting example of a first order bilinear invariant differential pairing in projective geometry, consider pairings
$$\odot^{2}T\mathcal{M}\times\Omega^{1}(\mathcal{M})\rightarrow \mathcal{O},$$
where $\odot^{2}$ denotes the second symmetric product and $\mathcal{O}$ is the sheaf of holomorphic (or smooth) functions. The transformation rule for $V^{ab}\in\odot^{2}T\mathcal{M}$ under change of connection is given by $\hat{\nabla}_{c}V^{ab}=\nabla_{c}V^{ab}+2\Upsilon_{c}V^{ab}+2\Upsilon_{d}V^{d(a}\delta_{c}{}^{b)}$. This implies
\begin{eqnarray*}
V^{ab}\hat{\nabla}_{(a}\omega_{b)}&=&V^{ab}\nabla_{(a}\omega_{b)}-2V^{ab}\Upsilon_{a}\omega_{b}\quad\text{and}\\
\omega_{b}\hat{\nabla}_{a}V^{ab}&=&\omega_{b}\nabla_{a}V^{ab}+(n+3)\omega_{b}\Upsilon_{a}V^{ab},
\end{eqnarray*}
where~$n=\mathrm{dim}(\mathcal{M})$. Therefore the pairing
$$(n+3)V^{ab}\nabla_{(a}\omega_{b)}+2\omega_{b}\nabla_{a}V^{ab}$$
does not depend upon a choice of connection within the projective class.
\par
\vspace{0.2cm}
It is natural to ask the question of whether these are all first order bilinear invariant differential pairings between those bundles and whether one can classify pairings between arbitrary bundles in general. This paper takes a first step towards a classification of those pairings and is divided into three parts:
\par
In the first part we construct first order bilinear invariant differential pairings on all manifolds with an AHS structure following closely and using the strategy and techniques developed in~\cite{css}. In the flat homogeneous (complex) case all manifolds with an AHS structure are of the form $G/P$, where $G$ is a complex semisimple Lie group and $P$ is a maximal parabolic subgroup. In this case the construction yields a generic classification of all first order bilinear invariant differential pairings between irreducible homogeneous bundles except for certain totally degenerate cases.
The first order case is made particularly easy by the fact that all first order linear invariant differential operators between so called irreducible associated bundles
are known and can be described as in ~\cite{css}, p.~70, corollary~7.2. Moreover there is no curvature involved for first order pairings, so the formulae for the general case do not differ from the ones in the model case of homogeneous spaces~$G/P$.
\par
In the second part we develop a strategy to construct higher order invariant differential pairings for homogeneous spaces $G/P$ with an AHS structure. The formal setup resembles the procedure in ~\cite{bceg} and the basic idea that lies behind the construction of the pairings in the second part comes from the Jantzen- Zuckermann translation principle and involves an argument about central character. 
\par
The third part of the paper is dedicated exclusively to the flat model case $G/P=\mathbb{CP}_{n}$ and we classify (with a minor restriction on the representation) all bilinear invariant differential pairings for non excluded weights, i.e.~weights that do not induce invariant operators emanating from the bundles in question (the notion of {\bf weight} will be made precise in~3.5.).
The reason for working on $\mathbb{CP}_{n}$ is that all linear invariant differential operators between irreducible homogeneous bundles are standard (i.e. they correspond to individual arrows in the BGG resolution, see ~\cite{bgg} and ~\cite{l}) and can be described as in ~~\cite{css}, p.~65, corollary~5.3 and p.~68, theorem~6.5. In fact if $P=B$ is a Borel subgroup, then
a classification of all linear invariant operators has been given dually in terms on Verma modules in~~\cite{bgg}. If $P$ is a general parabolic, then this is still an open problem, but in certain cases, like for $\mathbb{CP}_{n}$, a classification is given in ~\cite{bc1} and~~\cite{bc2}.

\section{Conventions}
\subsection{Composition series}
We will write composition series with the help of $+$ signs as explained in ~\cite{beg}, p.~1193 and ~\cite{es}, p.~434, so a short exact sequence 
$$0\rightarrow \mathbb{A}\rightarrow \mathbb{B}\rightarrow \mathbb{C}\rightarrow 0$$
of modules is equivalent to writing a composition series
$$\mathbb{B}=\mathbb{C}+\mathbb{A}.$$
This notation has the advantage that one can conveniently write down the subquotients of any filtration, so a composition series
$$\mathbb{B}=\mathbb{A}_{0}+\mathbb{A}_{1}+...+\mathbb{A}_{N}$$
denotes a filtration
$$\mathbb{A}_{N}=A_{N}\subseteq A_{N-1}\subseteq...\subseteq A_{0}=\mathbb{B},$$
with $\mathbb{A}_{i}=A_{i}/A_{i+1}$.
It can be noted that every composition series $\mathbb{B}=\mathbb{A}_{0}+\mathbb{A}_{1}+...+\mathbb{A}_{N}$ has a projection $\mathbb{B}\rightarrow \mathbb{A}_{0}$ and injections $\mathbb{A}_{j}+..+\mathbb{A}_{N}\rightarrow \mathbb{B}$, for $j=0,...,N$.

\subsection{Dynkin diagrams}
To denote a representation $(\mathbb{V},\rho)$ of a simple Lie algebra $\g$ or a parabolic subalgebra $\p\subset\g$, we write down the coefficient $B(\lambda,\alpha_{j}^{\vee})$ over the $j$-th node in the Dynkin diagram for $\g$, with $\lambda$ being the highest weight of the dual representation $(\mathbb{V}^{*},\rho^{*})$, $B(.,.)$ the Killing form and the $\alpha_{j}^{\vee}$ are the co-roots of the simple roots $\alpha_{0},...,\alpha_{n-1}$  of~$\g$. The details for this construction and the reason for this slightly odd notation is explained in ~\cite{be}, p.~22,23. To indicate that a representation of a parabolic subalgebra $\p$ (that corresponds to a subset 
 $\mathcal{S}_{\p}$ of the simple roots of the Lie algebra $\g$) is being referred to, we cross through all nodes in the Dynkin diagram for $\g$  that lie in~$\mathcal{S}\backslash\mathcal{S}_{\p}$.
 \par
If $(\mathbb{V},\rho)$ is a finite dimensional and irreducible representation of $\g$, then all the numbers above the nodes have to be non-negative integers. If  $(\mathbb{V},\rho)$ is a finite dimensional and irreducible representation of $\p$, then the numbers over all nodes have to be integers and the numbers over the uncrossed nodes have to be non-negative (these correspond to the irreducible representation of the semisimple part $\g_{0}^{S}$ of~$\p$). To be more precise, for a representation of $\p$ it would be sufficient to have non-negative integers over all the uncrossed nodes, but we will only deal with representations of $\p$ that lift to representations of  $P$, so we have to demand that the coefficients over the crossed through nodes are integers as well (see ~\cite{be}, p.~23, remark~3.1.6.).

\subsubsection{Example}
The bundle
$$\mathcal{O}(k)=\xo{k}{0}\;...\;\oo{0}{0}$$
denotes the $k$-th tensor power of the hyperplane section bundle on $\mathbb{CP}_{n}$ and is induced by the one dimensional representation of $\p$ on which $\h$, the Cartan subalgebra of $\g$, acts as
$$\left(\begin{array}{cc}
a & 0\\
0 & *
\end{array}\right).z=-kaz.$$

\par
The Dynkin diagram for an irreducible finite dimensional representation $\mathbb{E}$ of $\p$ will also be used to denote the corresponding {\bf generalized}
{\bf Verma} {\bf module}
$$M_{\p}(\mathbb{E})=\mathfrak{U}(\g)\otimes_{\mathfrak{U}(\p)}\mathbb{E}^{*},$$
where $\mathfrak{U}(\mathfrak{a})$ is the universal enveloping algebra for a Lie  algebra~$\mathfrak{a}$. More information about generalized Verma modules can be found in 
~\cite{l}, p.~500 and its correlation to irreducible homogeneous vector bundles on $G/P$ is explained in ~\cite{be}, p.~164. If $\mathbb{E}$ is just finite dimensional, then $M_{\p}(\mathbb{E})$ is called {\bf induced} {\bf module}.

\section{AHS structures}

\subsection{$|1|$-graded Lie algebras}
The basic ingredient of our construction is a complex semisimple Lie group $G$ with a $|1|$-graded Lie algebra
$$\g=\g_{-1}\oplus\g_{0}\oplus\g_{1},$$
where $[\g_{i},\g_{j}]\subset \g_{i+j}$. We will write $\p=\g_{0}\oplus\g_{1}$ and note that $\g_{\pm1}$ are commutative and dual with respect to the Killing form. Moreover $\g_{0}$ is reductive and has a semisimple part $\g_{0}^{S}=[\g_{0},\g_{0}]$ and a one dimensional centre, that is spanned by a grading element $E$, such that the decomposition of
$\g$ corresponds to a decomposition  into eigenspaces for the adjoint action of $E$, i.e.~$[E,X]=jX$ if and only if~$X\in\g_{j}$. More information and proofs about graded Lie algebras are to be found in~~\cite{cs}.
\par
For consistency reasons we will state all results in the holomorphic category but the construction carries over into the smooth category with minor modifications, see~~\cite{css}, p.~66).

\subsubsection{Example}
Complex projective $n$-space $\mathbb{CP}_{n}$ can be realized as $G/P$, with $G=\mathrm{SL}_{n+1}\C$ and Lie algebra $\g=\mathfrak{sl}_{n+1}\C$ with 
a $|1|$-grading
$$\g_{-1}=\left(\begin{array}{cc}
 0& 0\\
*&0\end{array}\right),\quad 
\g_{0}=\left(\begin{array}{cc}
 *& 0\\
0&*\end{array}\right),\quad
\g_{1}=\left(\begin{array}{cc}
 0& *\\
0&0\end{array}\right)\;\text{and}\;
\g_{0}^{S}=\left(\begin{array}{cc}
 0& 0\\
0&*\end{array}\right),\quad$$
where the sizes of the blocks  are $1\times1$, $1\times n$, $n\times 1$ and $n\times n$ from top left to bottom right. The grading element $E$ is given by
$$E=\left(\begin{array}{cc}
\frac{n}{n+1}& 0\\
0 & -\frac{1}{n+1}I_{n}
\end{array}\right),$$
where $I_{n}$ denotes the $n\times n$ standard matrix.

\subsection{Cartan connection}
The second ingredient is a {\bf Cartan} {\bf Geometry} $(\mathcal{G},\eta ,P,(\mathfrak{g},\mathfrak{p}))$ on  $\mathcal{M}$ consisting of
\begin{enumerate}
\item
a manifold $\mathcal{M}$,
\item
a principal right $P$ bundle $\mathcal{G}$ over $\mathcal{M}$ and
\item
a $\mathfrak{g}$-valued $1$-form $\eta$ on $\mathcal{G}$ satisfying the following conditions:
\begin{enumerate}
\item
the map $\eta _{g}:T_{g}\mathcal{G}\longrightarrow \mathfrak{g}$ is a linear isomorphism for every $g\in \mathcal{G}$,
\item
$R_{p}^{*}\eta =\mathrm{Ad}(p^{-1})\circ\eta $ for all $p\in P$, where $R_{p}$ denotes the natural right action of an element $p\in P$ in the structure group,  and
\item
$\eta (\zeta_{X})=X$ for all $X\in \mathfrak{p}$, where $\zeta_{X}$ is the (vertical) fundamental vector field on $\mathcal{G}$ associated to~$X\in\p$.
\end{enumerate}
\end{enumerate}
Examples of Cartan Geometries are given by the {\bf flat} {\bf models} $\mathcal{M}=G/P$ with the {\bf Maurer} {\bf Cartan} {\bf form} $\eta_{g}=(L_{g^{-1}})_{*}$, where $L_{g}$ denotes left multiplication. We will deal with a very specific example of such a flat model in the third part of the paper.

\subsection{Associated bundles} 
For every finite dimensional representation $\rho:P\rightarrow \mathrm{Aut}(\mathbb{V})$, we can define a corresponding {\bf associated} vector bundle
$$V=(V,\rho)=\mathcal{G}\times_{P}\mathbb{V}=(\mathcal{G}\times\mathbb{V})/\sim,$$
where
$$(gp,v)\sim(g,\rho(p)v)\;\forall\;g\in \mathcal{G},p\in P\;\text{and}\;v\in\mathbb{V}.$$
These bundles are exactly the {\bf homogeneous} {\bf bundles} in the flat case $G/P$ (in general this procedure is functorial and defines a {\bf natural} {\bf vector} {\bf bundle} in the sense of~~\cite{kms}).
Sections of this bundle can be identified with maps
$$s:\mathcal{G}\rightarrow\mathbb{V}\,\text{s.t.}\;s(gp)=\rho(p^{-1})s(g),$$
for all $g\in\mathcal{G}$ and~$p\in P$. We can differentiate this requirement to obtain
$$(\zeta_{X}s)(g)=-\rho(X)s(g)\;\forall\;X\in\p,g\in\mathcal{G}.$$
We will write $\Gamma(V)=\mathcal{O}(\mathcal{G},\mathbb{V})^{P}$ for the space of sections of~$V$. The tangent bundle $T\mathcal{M}$ and cotangent bundle $\Omega^{1}(\mathcal{M})$, for example, arise via the adjoint representation of $P$ on $\g/\p\cong\g_{-1}$ and its dual~$\g_{-1}^{*}\cong\g_{1}$. We will denote the bundles and sections of these bundles by the Dynkin diagram notation for the representation that induces them. In the case of $M=\mathbb{CP}_{n}$, for example,  we write
$$\xoo{1}{0}{0}\;...\;\oo{0}{1}=T\mathcal{M}\;\text{and}\;\xoo{-2}{1}{0}\;...\;\oo{0}{0}=\Omega^{1}(\mathcal{M}).$$

\subsection{The invariant differential}
The Cartan connection does not yield a connection on associated bundles, but we can still define the {\bf invariant} {\bf differential}
$$\nabla^{\eta}:\mathcal{O}(\mathcal{G},\mathbb{V})\rightarrow \mathcal{O}(\mathcal{G},\g_{-1}^{*}\otimes\mathbb{V})$$
with
$$\nabla^{\eta}s(g)(X)=\nabla^{\eta}_{X}s(g)=[\eta^{-1}(X)s](g),\;\forall\;X\in \g_{-1},g\in\mathcal{G}\;\text{and}\;s\in\mathcal{O}(\mathcal{G},\mathbb{V}).$$
It has to be noted that it does not take $P$-equivariant sections to $P$-equivariant sections.

\subsection{Representations of $\p$}
Since $\g_{1}$ is nilpotent, it acts trivially on any irreducible representation of $\p$ and so we can denote any irreducible representation of $\p$ by the highest weight of the corresponding representation of $\g_{0}^{S}$ and by specifying how $E$ acts. Since $E$ lies in the centre of $\g_{0}$, it will always act
by multiplication of a constant which we call, following ~\cite{cd}, {\bf geometric} {\bf weight} (or sometimes for brevity just weight).
The tangent and co-tangent space $\g_{-1}$ and $\g_{1}$ have geometric weight $-1$ and $1$ respectively.
\par
If $\h$, the Cartan subalgebra of $\g$, is chosen in such a way that all positive roots spaces lie in $\g_{0}\cup\g_{1}$ and so that $E\in\h$, then the Cartan subalgebra of $\g_{0}^{S}$ is given by $\h^{S}=\h\cap\g_{0}^{S}$ and we can denote representations of $\g_{0}$ by their highest weight in~$\h^{*}$. The representation is finite dimensional and irreducible if the restriction of this weight to $(\h^{S})^{*}$ is dominant integral.

\subsection{Jet bundles and invariant pairings}
For every complex (or smooth) manifold $\mathcal{M}$ and holomorphic vector bundle $V$ over $\mathcal{M}$, we denote by $J^{k}V$ the vector bundle over $\mathcal{M}$ of {\bf $k$-jets} of~$V$.
The fibre of $J^{k}V$ over each point $x\in \mathcal{M}$ is the quotient of the space of germs of sections of $V$ at $x$ by the subspace of germs of sections which vanish to order $k+1$ at~$x$. Linear differential operators $D:V\rightarrow W$ of order $k$ between two vector bundles are in one-to-one correspondence  with bundle homomorphisms $d:J^{k}V\rightarrow W$ (see ~\cite{s},~p.~183). This motivates one to define a {\bf bilinear} {\bf differential} {\bf pairing} between sections of bundles $V$ and $W$ to sections of a bundle $U$ by a homomorphism
$$d:J^{k}V\otimes J^{l}W\rightarrow U.$$
This pairing is of {\bf order $M$} if and only if 
\begin{enumerate}
\item
$k=l=M$,
\item
there is a subbundle $B$ of $J^{M}V\otimes J^{M}W$, so that there is a commutative diagram
$$\begin{array}{ccc}
J^{M}V\otimes J^{M}W &\\
\downarrow &\searrow d\\
 (J^{M}V\otimes J^{M}W)/B& \stackrel{\phi}{\rightarrow}& U
 \end{array}$$
and
\item
the map $\phi$ induces a formula that consist of terms in which derivatives of sections of $V$ are combined with derivatives of sections of $W$ in such a way that the total order is $M$ (i.e. a term may consist of a $k$-th derivative of a section of $V$ combined with a $(M-k)$-th derivative of a section of $W$, for~$k=0,...,M$).
\end{enumerate}
Note that there is always a canonical choice of $B$ as detailed in the Appendix. We will therefore write $J^{M}(V,W)=(J^{M}V\otimes J^{M}W)/B$ for the canonical choice of $B$.
This is not to be confused with the set of all $M$-jets of $V$ into $W$ as defined in ~\cite{kms}, p.~117, definition~12.2.
\par
If $\mathcal{M}=G/P$ is a homogeneous space, then a pairing is called {\bf invariant} (some authors use the term equivariant) if it commutes with the action of $G$ on sections of the involved homogeneous vector bundles, which is given by $(g.s)(h)=s(g^{-1}h)$ for all $g,h\in G$ and $s\in\Gamma(F)$ (see also lemma 1).
\par
In general, there is no commonly accepted notion of invariance for manifolds with an AHS (or more generally parabolic) structure (see ~\cite{sl}, p.~193, section 2). We will deal with this issue by taking a pragmatic point of view: First of all, every manifold with an AHS structure is equipped with a distinguished class of connections (Weyl connections), as detailed in ~\cite{css1},~p.~42 and~~\cite{css},~p.~54. A pairing is then called {\bf invariant}, if $\phi$ induces a formula that consists of terms involving an arbitrary connection from the distinguished equivalence class, but that as a whole does not depend on its choice.
This slightly delicate point will not play any role in this paper since we will only be constructing first order pairings on a general manifold with AHS structure, see lemma 1.

\par\vspace{0.3cm}
This is the general situation that we will be working in and we will treat the data described above as a given AHS structure on our manifold~$\mathcal{M}$.

\section{The first order case}
In the following, we will fix two finite dimensional irreducible representations $\tilde{\lambda}:\p\rightarrow\mathfrak{gl}(\mathbb{V})$ and $\tilde{\nu}:\p\rightarrow\mathfrak{gl}(\mathbb{W})$  that are
induced from irreducible representations of $\g_{0}^{S}$ with highest weights $\lambda$ and $\nu$ respectively and where the grading element $E$ acts as $\omega_{1}$ and
$\omega_{2}$ respectively. 
\par
\vspace{0.2cm}
In the homogeneous case $G/P$, the first {\bf jet} {\bf bundle} $J^{1}V$ associated to any homogeneous vector bundle is also homogeneous:  The fibre $J^{1}\mathbb{V}$ of $J^{1}V$ at the origin $P\in G/P$ consists of germs of sections of $V$ modulo those that vanish at $P$ of order at least 2. This vector space carries a representation of $P$ that makes $J^{1}V\cong G\times_{P}J^{1}\mathbb{V}$ a homogeneous bundle. In general, we can use this representation of $P$ on the \emph{vector space} $J^{1}\mathbb{V}=\mathbb{V}\oplus(\g_{1}\otimes \mathbb{V})$ to define an associated bundle that is exactly the first jet bundle $J^{1}V$ of $V$ as defined above, 
see~~\cite{css1},~p.~56. This construction ensures that the map
$$\mathcal{O}(\mathcal{G},\mathbb{V})^{P}\ni s\mapsto (s,\nabla^{\eta}s)\in\mathcal{O}(\mathcal{G},J^{1}\mathbb{V})^{P}$$
is well defined, i.e.~maps $P$ equivariant sections to $P$ equivariant sections, as shown in~~\cite{css1},~p.~56.
\par
There are two exact sequences associated to the first jet bundles of $V$ and $W$:
$$0\rightarrow \Omega^{1}\otimes V\rightarrow J^{1}V \rightarrow V\rightarrow 0$$
and 
$$0\rightarrow \Omega^{1}\otimes W\rightarrow J^{1}W\rightarrow W\rightarrow 0,$$
which are the {\bf jet} {\bf exact} {\bf sequences} as described in~~\cite{s},~p.~182. 
All these are associated bundles, so on the level of $\p$ representations we have
two filtered $\p$-modules
$$J^{1}\mathbb{V}=\mathbb{V}+\g_{1}\otimes \mathbb{V}$$
and
$$J^{1}\mathbb{W}=\mathbb{W}+\g_{1}\otimes \mathbb{W}.$$
Hence the tensor product has a filtration
$$J^{1}\mathbb{V}\otimes J^{1}\mathbb{W}=\mathbb{V}\otimes \mathbb{W}+\begin{array}{c}
\mathbb{V}\otimes\g_{1}\otimes \mathbb{W}\\
\oplus\\
\g_{1}\otimes \mathbb{V}\otimes \mathbb{W}
\end{array}
+\g_{1}\otimes \mathbb{V}\otimes\g_{1}\otimes \mathbb{W}.$$
The module
$$J^{1}(\mathbb{V},\mathbb{W})=J^{1}\mathbb{V}\otimes J^{1}\mathbb{W}/(\g_{1}\otimes \mathbb{V}\otimes\g_{1}\otimes \mathbb{W})$$
can, as a \emph{vector space}, be written as
$$J^{1}(\mathbb{V},\mathbb{W})=(\mathbb{V}\otimes\mathbb{W})\oplus(\mathbb{V}\otimes\g_{1}\otimes\mathbb{W})\oplus(\g_{1}\otimes\mathbb{V}\otimes\mathbb{W}).$$
The $\p$-module structure of $J^{1}(\mathbb{V},\mathbb{W})$ is, however, induced by the $\p$-module structures of $J^{1}\mathbb{V}$ and $J^{1}\mathbb{W}$ that induce the associated bundles $J^{1}V$ and $J^{1}W$. It is defined in such a way that the mapping
$$\mathcal{O}(\mathcal{G},\mathbb{V})^{P}\otimes \mathcal{O}(\mathcal{G},\mathbb{W})^{P}\ni(s,t)\mapsto
(s\otimes t,s\otimes \nabla^{\eta}t,\nabla^{\eta}s\otimes t)\in\mathcal{O}(\mathcal{G},J^{1}(\mathbb{V},\mathbb{W}))^{P}$$
is well defined. More precisely, it can be checked directly, following the strategy in ~\cite{css1},~p.~56, that the action is given by


\begin{eqnarray*}
&&j^{1}(\tilde{\lambda},\tilde{\nu})(Z)[(v,\varphi)\otimes(w,\phi)]\\
&=&\left(\begin{array}{c}
\tilde{\lambda}(Z)v\otimes w+v\otimes\tilde{\nu}(Z)w\\
\tilde{\lambda}(Z)v\otimes \phi+v\otimes(\tilde{\nu}(Z)\circ\phi-\phi\circ\mathrm{ad}_{\g_{-1}}(Z)+\tilde{\nu}(\mathrm{ad}_{\p}(Z)(...))w)\\
\varphi\otimes\tilde{\nu}(Z)w+(\tilde{\lambda}(Z)\circ \varphi-\varphi\circ\mathrm{ad}_{\g_{-1}}(Z)+\tilde{\lambda}(\mathrm{ad}_{\p}(Z)(...))v)\otimes w
\end{array}\right),
\end{eqnarray*}
where $(v,\varphi)\in J^{1}\mathbb{V}$, $(w,\phi)\in J^{1}\mathbb{W}$, the square brackets denote projection onto $J^{1}(\mathbb{V},\mathbb{W})$ and the expressions $\mathrm{ad}_{\mathfrak{a}}$ for a subalgebra $\mathfrak{a}$ of $\g$ stand for the usual adjoint representation followed by the projection onto~$\mathfrak{a}$.
Indeed we must have (see~3.3)
$$\eta^{-1}(Z).(s\otimes t,s\otimes \nabla^{\eta}_{X}t,\nabla^{\eta}_{Y}s\otimes t)=-j^{1}(\tilde{\lambda},\tilde{\nu})(Z).(s\otimes t,s\otimes \nabla^{\eta}_{X}t,\nabla^{\eta}_{Y}s\otimes t),$$
for all $Z\in\p$ and~$X,Y\in \g_{-1}$.
Since we can compute the left hand side, the right hand side and therefore the representation $j^{1}(\tilde{\lambda},\tilde{\nu})$ is uniquely defined. The reason for this setup is given in the following lemma.

\begin{lemma}
First order bilinear invariant differential pairings
$$\Gamma(V)\times\Gamma(W)\rightarrow \Gamma(E)$$
in the flat homogeneous case $G/P$ are in one-to-one correspondence with $\p$-module homomorphisms
$$J^{1}(\mathbb{V},\mathbb{W})\rightarrow \mathbb{E}.$$
In the general AHS case, these homomorphisms yield first order bilinear differential pairings, which are invariant in the sense that they produce formulae which do not depend on 
a particular choice of connection within the distinguished class.
\end{lemma}\begin{proof}
A first order linear differential operator $\phi:\Gamma(V)\rightarrow\Gamma(F)$ corresponds to a homomorphism  $J^{1}V\rightarrow F$ (by definition of the jet bundles).
In the flat homogeneous case $G/P$, an operator is called invariant if it commutes with the action of $G$, hence those operators are uniquely determined by their action at the identity coset $P$, where they induce $\p$-module homomorphisms
$$J^{1}\mathbb{V}\rightarrow \mathbb{F}.$$
Conversely, every such $\p$-module homomorphism induces a first order invariant linear differential operator. 
Analogous reasoning shows that first order invariant bilinear differential pairings as defined in 3.6 are in one-to-one correspondence with $\p$-module homomorphisms $J^{1}\mathbb{V}\otimes J^{1}\mathbb{W}\rightarrow \mathbb{E}$ that factor through an appropriate subbundle.
This subbundle is $\g_{1}\otimes\mathbb{V}\otimes\g_{1}\otimes\mathbb{W}$, because expressions derived from this subbundle correspond to terms with two derivatives and they have an incorrect geometric weight. To be more precise, the irreducible components of $\g_{1}\otimes\mathbb{V}\otimes\g_{1}\otimes\mathbb{W}$ have geometric weight $\omega_{1}+\omega_{2}+2$, whereas the irreducible components of $\g_{1}\otimes\mathbb{V}\otimes\mathbb{W}$ have geometric weight $\omega_{1}+\omega_{2}+1$, due to the fact that $\g_{1}$ has geometric weight 1. This implies that in formulae we will only be allowed to use terms that have derivatives in either
sections of $V$ or sections of $W$ but not in both. Algebraically this means that we can factor out $\g_{1}\otimes\mathbb{V}\otimes\g_{1}\otimes\mathbb{W}$  of~$J^{1}\mathbb{V}\otimes J^{1}\mathbb{W}$.
\par
In the general AHS case, as mentioned above, the $\p$-modules structure of $J^{1}(\mathbb{V},\mathbb{W})$ ensures that the mapping
$$\mathcal{O}(\mathcal{G},\mathbb{V})^{P}\otimes \mathcal{O}(\mathcal{G},\mathbb{W})^{P}\ni(s,t)\mapsto
(s\otimes t,s\otimes \nabla^{\eta}t,\nabla^{\eta}s\otimes t)\in\mathcal{O}(\mathcal{G},J^{1}(\mathbb{V},\mathbb{W}))^{P}$$
is well defined. Following ~\cite{css}, p.~54, we note that the Cartan connection is uniquely defined by the AHS structure, so the differential pairings that we obtain from $\p$-module homomorphisms $J^{1}(\mathbb{V},\mathbb{W})\rightarrow \mathbb{E}$ combined with this map do not depend upon a specific choice of connection within the distinguished class (see ~\cite{sl}, p.~194, 197). Pairings that arise via this construction are {\bf strongly invariant} in the sense of ~\cite {css2}, p.~102. 
More information about canonical Cartan connections associated to AHS structures can be found in~~\cite{css1}.
\end{proof}\par
\vspace{0.3cm}

Looking at the exact sequence of $\p$-modules
$$0\rightarrow\begin{array}{c}
\g_{1}\otimes \mathbb{V}\otimes \mathbb{W}\\
\oplus\\
\mathbb{V}\otimes\g_{1}\otimes \mathbb{W}
\end{array}\rightarrow J^{1}(\mathbb{V},\mathbb{W})\rightarrow  \mathbb{V}\otimes\mathbb{W}\rightarrow 0,$$
it is clear that a $\p$-module homomorphism $J^{1}(\mathbb{V},\mathbb{W})\rightarrow\mathbb{E}$ onto an irreducible $\p$-module $\mathbb{E}$ induces a $\g_{0}^{S}$-homomorphism
$$\begin{array}{c}
\g_{1}\otimes \mathbb{V}\otimes \mathbb{W}\\
\oplus\\
\mathbb{V}\otimes\g_{1}\otimes \mathbb{W}
\end{array}\stackrel{\pi}{\rightarrow} \mathbb{E}$$
and so the only candidates for $\mathbb{E}$ are the irreducible components of
$\g_{1}\otimes\mathbb{V}\otimes\mathbb{W}$ viewed as $\g_{0}^{S}$-modules.
However, not every projection $\pi$ is a $\p$-module homomorphism. In order to determine which $\pi$ are allowed, it can be noted that the action of $\g_{0}$ on $J^{1}(\mathbb{V},\mathbb{W})$ is just the tensorial one, so $J^{1}(\mathbb{V},\mathbb{W})$ can be split as a $\g_{0}$-module. But $\g_{1}$ does not act trivially as on any irreducible $\p$-module, so in order to check that a specific projection is indeed a $\p$-module homomorphism and not just
a $\g_{0}^{S}$-module homomorphism the image of the action of $\g_{1}$, when acting in $J^{1}(\mathbb{V},\mathbb{W})$, has to vanish under~$\pi$.  
On the other hand this is obviously sufficient for $\pi$ to be a $\p$-module
homomorphism. We therefore compute for $Z\in \g_{1}$:
$$j^{1}(\tilde{\lambda},\tilde{\nu})(Z)[(v,\varphi)\otimes (w,\phi)]=(0,v\otimes\tilde{\nu}(\mathrm{ad}_{\p}(Z)(...))w,\tilde{\lambda}(\mathrm{ad}_{\p}(Z)(...))v\otimes w)$$
and call this term {\bf obstruction} {\bf term}.
\par
\vspace{0.2cm} 

In order to obtain an explicit formula for the obstruction term we will rewrite this expression: 
The action of an element $Z\in\g_{1}$ in $J^{1}(\mathbb{V},\mathbb{W})$ can be interpreted as a map
$$Z:\mathbb{V}\otimes\mathbb{W}\rightarrow 
\begin{array}{c}
\mathbb{V}\otimes\g_{1}\otimes\mathbb{W}\\
\oplus\\
\g_{1}\otimes\mathbb{V}\otimes\mathbb{W}\\
\end{array}
\cong
\begin{array}{c}
\mathrm{Hom}(\g_{-1},\mathbb{V}\otimes\mathbb{W})\\
\oplus\\
\mathrm{Hom}(\g_{-1},\mathbb{V}\otimes\mathbb{W})
\end{array},$$
with
$$Z(v\otimes w)\left(\begin{array}{c}
X\\
Y
\end{array}\right)=\left(
\begin{array}{c}
v\otimes \tilde{\nu}([Z,X])w\\
\tilde{\lambda}([Z,Y])v\otimes w
\end{array}\right).$$
For the following it is convenient to normalize the Killing form $B(.,.)$ to a form $(.,.)$ with $(E,E)=1$, where $E$ is the grading element.
Having done this, we introduce dual basis $\{\eta_{\alpha}\}$ and $\{\xi_{\alpha}\}$ of $\g_{1}$ and $\g_{-1}$ respectively with respect to this form. This yields
$$X=\sum_{\alpha}\eta_{\alpha}(X)\xi_{\alpha}$$
for every $X\in\g_{-1}$. Writing down the obstruction term as a mapping
$$\Phi:\g_{1}\otimes\mathbb{V}\otimes\mathbb{W}\rightarrow
\begin{array}{c}
\mathbb{V}\otimes\g_{1}\otimes\mathbb{W}\\
\oplus\\
\g_{1}\otimes\mathbb{V}\otimes\mathbb{W}\\
\end{array},\quad Z\otimes v\otimes w\mapsto Z(v\otimes w), $$
with
$$\Phi(Z\otimes v\otimes w)=\left(
\begin{array}{c}
v\otimes \sum_{\alpha}\eta_{\alpha}\otimes\tilde{\nu}([Z,\xi_{\alpha}])w\\
\sum_{\alpha}\eta_{\alpha}\otimes\tilde{\lambda}([Z,\xi_{\alpha}])v\otimes w
\end{array}\right),$$
allows one to use the casimir operator to turn this into an easier expression given by the following lemma.

\begin{lemma}
A projection 
$$\pi:\begin{array}{c}
\mathbb{V}\otimes\g_{1}\otimes\mathbb{W}\\
\oplus\\
\g_{1}\otimes\mathbb{V}\otimes\mathbb{W}
\end{array}
\rightarrow \mathbb{E}$$
onto an irreducible component $\mathbb{E}$ of the $\g_{0}^{S}$ tensor product $\mathbb{V}\otimes\g_{1}\otimes\mathbb{W}$ is a $\p$-module homomorphism if and only if~$\pi\circ\Phi=0$. The mapping $\Phi$ can be written as
$$\Phi(Z\otimes v \otimes w)=\left(
\begin{array}{c}
v\otimes\sum_{\sigma}(\omega_{2}-c_{\nu\sigma})\pi_{\nu\sigma}(Z\otimes w)\\
\sum_{\tau}(\omega_{1}-c_{\lambda\tau})\pi_{\lambda\tau}(Z\otimes v)\otimes w
\end{array}\right),$$
where $\tau$ and $\sigma$ range over the highest weights of the irreducible components of $\g_{1}\otimes\mathbb{V}$ and $\g_{1}\otimes\mathbb{W}$ respectively
and $\pi_{\lambda\tau}$, $\pi_{\nu\sigma}$ denote the corresponding projections. The constants $c_{\gamma\delta}$ are defined by
$$c_{\gamma\delta}= -\frac{1}{2}[(\delta,\delta+2\rho)-(\gamma,\gamma+2\rho)-(\alpha,\alpha+2\rho)],$$
where 
$\alpha$ is the highest weight of $\g_{1}$ and 
$$\rho=\rho_{\g_{0}^{S}}=\frac{1}{2}\sum_{\beta\in\Delta^{+}(\g_{0}^{S})}\beta,$$
with $\Delta^{+}(\g_{0}^{S})$ denoting the set of positive roots of $\g_{0}^{S}$.
\end{lemma}
\begin{proof}
The first calculation in this direction in the conformal case was done in ~\cite{f} and the general case is proved in~~\cite{css},~p.~63,~lemma~4.3.
\end{proof}
\par
\vspace{0.4cm}
The decompositions $\g_{1}\otimes\mathbb{V}$ and $\g_{1}\otimes\mathbb{W}$ do not have multiplicities  (see~~\cite{css},~p.~58,59), so let us write
$$\g_{1}\otimes\mathbb{V}=\mathbb{V}(\tau_{1})\oplus...\oplus\mathbb{V}(\tau_{r})$$
and
$$\g_{1}\otimes\mathbb{W}=\mathbb{V}(\sigma_{1})\oplus...\oplus\mathbb{W}(\sigma_{s}),$$
where the greek letter in the brackets denotes the highest weight of the module. These weights all correspond to the action of $\g_{0}^{S}$, the geometric weight of $\mathbb{V}(\tau_{i})$ is $\omega_{1}+1$ for  every $i$, since it lies in the tensor product $\g_{1}\otimes\mathbb{V}$ and $\g_{1}$ has geometric weight~1.
Analogously, the geometric weight of all $\mathbb{W}(\sigma_{j})$ is~$\omega_{2}+1$.
\par
If $\mathbb{E}$ is one of the irreducible components of $\g_{1}\otimes\mathbb{V}\otimes\mathbb{W}$ of highest weight $\mu$, then we denote by
$\pi_{\tau\mu}^{i}$ the projection $\mathbb{V}({\tau})\otimes\mathbb{W}\rightarrow \mathbb{E}^{(i)}$ into the $i$-th copy of $\mathbb{E}$ in
the decomposition. $\pi_{\sigma\mu}^{j}$ is defined analogously as the projection into the $j$-th copy of $\mathbb{E}$ in $\mathbb{V}\otimes\mathbb{W}({\sigma})$.
Every projection 
$$\pi:
\begin{array}{c}
\mathbb{V}\otimes\g_{1}\otimes\mathbb{W}\\
\oplus\\
\g_{1}\otimes\mathbb{V}\otimes\mathbb{W}
\end{array}
\rightarrow \mathbb{E}$$
can be written as
\begin{eqnarray*}
\pi \left(
\begin{array}{c}
v_{1}\otimes Z_{1}\otimes w_{1}\\
Z_{2}\otimes v_{2}\otimes w_{2} \\
\end{array}\right)
& = &
\sum_{\tau}\sum_{i}a_{\tau,i}\pi^{i}_{\tau\mu}\left(\pi_{\lambda\tau}(Z_{1}\otimes v_{1})\otimes w_{2}\right)\\
&   & +
\sum_{\sigma}\sum_{j}b_{\sigma,j}\pi_{\sigma\mu}^{j}\left(v_{2}\otimes \pi_{\nu\sigma}(Z_{2}\otimes w_{2})\right),
\end{eqnarray*}
for some constants $a_{\tau,i}$ and $b_{\sigma,j}$.
In order for a projection $\pi$ to be a $\p$-homomorphism, $\pi\circ\Phi(Z\otimes v\otimes w))=0$ has to hold for all
$Z\in\g_{1}$, $v\in\mathbb{V}$ and~$w\in\mathbb{W}$. This reads
\begin{eqnarray*}
\pi\circ\Phi (Z\otimes v\otimes w)& = &
\sum_{\tau}\sum_{i}a_{\tau,i}(\omega_{1}-c_{\lambda\tau})\pi^{i}_{\tau\mu}(\pi_{\lambda\tau}(Z\otimes v)\otimes w)\\
&   & +\sum_{\sigma}\sum_{j}b_{\sigma,j}(\omega_{2}-c_{\nu\sigma})\pi^{j}_{\sigma\mu}(v\otimes\pi_{\nu\sigma}(Z\otimes w))\\
& = & 0.
\end{eqnarray*}
Let $k$ denote the number of copies of $\mathbb{E}$ in $\g_{1}\otimes \mathbb{V}\otimes\mathbb{W}$, then there are 
$2k$ unknowns and $k$ equations. Since $Z$, $v$ and $w$ are to be arbitrary and all $\pi_{\tau\mu}^{i}(\pi_{\lambda\tau}(Z\otimes v)\otimes w)$ lie in different
copies of $\mathbb{E}$, we can think of those elements as constituting a basis $\{e_{i}\}$ of~$\oplus^{k}\mathbb{E}$. The same is true for the different
$\pi_{\sigma\mu}^{j}(v\otimes\pi_{\nu\sigma}(Z\otimes w))$, which constitute a different basis~$\{f_{j}\}$. Hence there is a linear isomorphism 
$f_{j}=\sum_{i}A_{ij}e_{i}$ connecting those two basis and we obtain $k$ equations
$$a_{i}(\omega_{1}-c_{\lambda\tau(i)})+\sum_{j}b_{j}(\omega_{2}-c_{\nu\sigma(j)})A_{ij}=0,\; i=1,..,k,$$
where $\tau(i)$ (resp. $\sigma(j)$) denotes the representation corresponding to the index $i$ (resp. $j$), i.e. the $i$-th (resp. $j$-th) copy of $\mathbb{E}$ lies in $\mathbb{V}(\tau(i))\otimes\mathbb{W}$ (resp. in~$\mathbb{V}\otimes\mathbb{W}(\sigma(j))$). If $\omega_{1}$ does not equal one of the excluded weights 
 $c_{\lambda\tau}$, then the constants $a_{i}$ are uniquely determined by the~$b_{j}$'s. 
This yields a $k$-parameter family of invariant differential pairings if the geometric weight $\omega_{1}$ is  not 
excluded. An excluded geometric weight $\omega_{1}$ corresponds to an invariant differential operator 
$$\Gamma(V)\rightarrow\Gamma(V(\tau)),$$
where $V(\tau)$ is induced from the  representation $\mathbb{V}(\tau)$ with~$\omega_{1}=c_{\lambda\tau}$. The roles of the $a_{i}$ and $b_{j}$ can, of course, be interchanged, so that we can alternatively exclude geometric weights $\omega_{2}$, which correspond to first order invariant differential operators $\Gamma(W)\rightarrow\Gamma(W(\sigma))$. Thus we have proved:


\begin{theorem}[Main Result 1]
Let $\mathbb{V}$ and $\mathbb{W}$ be two irreducible $\p$-modules with decompositions
$$\g_{1}\otimes\mathbb{V}=\mathbb{V}(\tau_{1})\oplus...\oplus\mathbb{V}(\tau_{r})$$
and
$$\g_{1}\otimes\mathbb{W}=\mathbb{W}(\sigma_{1})\oplus...\oplus\mathbb{W}(\sigma_{s}).$$
If $\omega_{1}\not\in\{c_{\lambda\tau_{1}},...,c_{\lambda\tau_{r}}\}$ or $\omega_{2}\not\in\{c_{\nu\sigma_{1}},...,c_{\nu\sigma_{s}}\}$, then
there exists a $k$-parameter family of first order invariant bilinear differential pairings
$$\Gamma(V)\times\Gamma(W)\rightarrow \Gamma(E),$$
where $k$ is the number of copies of $\mathbb{E}$ in~$\g_{1}\otimes\mathbb{V}\otimes\mathbb{W}$. These are the only possible first order invariant bilinear differential pairings between section of $V$ and $W$ onto an irreducible bundle in the flat homogeneous case~$G/P$.
\end{theorem}

\subsubsection{Remark}
In fact only those weights $c_{\lambda\tau}$ (resp. $c_{\nu\sigma}$) have to be excluded for which $\mathbb{E}\subset\mathbb{V}(\tau)\otimes\mathbb{W}$ 
(resp.~$\mathbb{E}\subset \mathbb{V}\otimes\mathbb{W}(\sigma)$).

\begin{corollary}
The situation is considerably simplified if there is only one copy of $\mathbb{E}$ in~$\g_{1}\otimes\mathbb{V}\otimes\mathbb{W}$. 
Then we can choose $a=(\omega_{2}-c_{\nu\sigma})$ and $b=-(\omega_{1}-c_{\lambda\tau})$ if we normalize the projections correctly. Every multiple of this pairing is obviously invariant as well. It also shows what happens if weights are excluded:
\begin{enumerate}
\item If $\omega_{1}=c_{\lambda\tau}$, then we must take $b=0$ and $a$ is arbitrary. This corresponds to an invariant first order linear differential operator $\Gamma(V)\rightarrow \Gamma(V(\tau))$ combined with a projection~$\Gamma(V(\tau))\times\Gamma(W)\rightarrow \Gamma(E)$.
\item If $\omega_{2}=c_{\nu\sigma}$, then there is a first order linear invariant differential operator $\Gamma(W)\rightarrow\Gamma(W(\sigma))$ that can be combined with a projection $\Gamma(W(\sigma))\times\Gamma(V)\rightarrow\Gamma(E)$, i.e. we must take $a=0$ and $b$ is arbitrary. 
\item If both weights are excluded, then the statement of the main theorem is not true anymore. We obtain two independent pairings corresponding to the two invariant differential operators and the projections mentioned above.
\end{enumerate}
\end{corollary}

\subsubsection{Examples}
Looking at pairings on $\mathbb{CP}_{n}$ involves $\g=\mathfrak{sl}_{n+1}\C$ and~$\g_{0}^{S}\cong\mathfrak{sl}_{n}\C$.

\begin{enumerate}
\item
Let
$$V=\;\xoo{w}{0}{0}\;...\;\oo{0}{0}=\mathcal{O}(w)\quad\text{and}\quad W=\;\xoo{1+v}{0}{0}\;...\;\oo{0}{1},$$
denote weighted functions and weighted vector fields on~$\mathbb{CP}_{n}$. The tensor product decomposes as
$$\g_{1}\otimes\mathbb{V}\otimes\mathbb{W}=\;\oo{1}{0}\;...\;\oo{0}{1}\;\oplus\;\oo{0}{0}\;...\;\oo{0}{0}$$
and the weights are given by $\omega_{1}=-w\frac{n}{n+1}$, $\omega_{2}=-\frac{nv+n+1}{n+1}$. Taking $\mu=0$ yields $c_{\lambda\tau} =0$ and~$c_{\nu\mu}=n-1$. This corresponds to the invariant pairing (where we have multiplied everything by~$-\frac{n+1}{n}$):
$$(n+v+1)X^{a}\nabla_{a}f-w(\nabla_{a}X^{a})f.$$
\item
Quite similarly we obtain an invariant paring
$$\Omega^{1}(v)\times\mathcal{O}(w)\ni (\sigma_{b},f)\mapsto (v-2)\sigma_{(a}\nabla_{b)}f-w( \nabla_{(a}\sigma_{b)})f$$
from the fact that in this case $\omega_{2}-c_{\nu\sigma}=-\frac{n}{n+1}(v-2)$.
\item
A more sophisticated example can be obtained when we take
$$V=\;\xoo{1+v}{0}{0}\;...\;\oo{0}{1},\quad W=\underbrace{\;\xoo{w-(k+1)}{0}{0}\;...\;\ooo{0}{1}{0}\;...\;\oo{0}{0}}_{1\;\text{is in the}\;(k+1)\text{-th position}}$$
and
$$E=\;\underbrace{\xoo{v+w-(k+1)}{0}{0}\;...\;\ooo{0}{1}{0}\;...\;\oo{0}{0}}_{1\;\text{is in the}\;(k+1)\text{-th position}},$$
i.e.~we pair weighted vector fields with weighted $k$-forms to obtain weighted $k$-forms again. This time the multiplicity is two and indeed, for non-excluded weights, there is a two parameter family given by
$$X^{a}\nabla_{a}\omega_{bc...d}+\frac{n+v-w-vw+vk+1}{(n+v+1)(v+1)}(\nabla_{a}X^{a})\omega_{bc...d}-\frac{k+1}{v+1}(\nabla_{[a}X^{a})\omega_{bc...d]}$$
and
$$X^{a}\nabla_{[a}\omega_{bc...d]}+\frac{(n-k)w}{(n+v+1)(v+1)(k+1)}(\nabla_{a}X^{a})\omega_{bc...d}-\frac{w}{v+1}(\nabla_{[a}X^{a})\omega_{bc...d]}.$$
The denominators can only be zero, if the weights are excluded, because $\omega_{1}-c_{\lambda\tau_{1}}=-\frac{n}{n+1}(n+v+1)$ and
$\omega_{1}-c_{\lambda\tau_{2}}=-\frac{n}{n+1}(v+1)$. If one of these is zero, then the corresponding operator
$X^{a}\mapsto \nabla_{a}X^{a}$ or $X^{a}\mapsto\nabla_{b}X^{a}-\frac{1}{n}\nabla_{c}X^{c}\delta_{b}{}^{a}$ is projectively invariant.
If we take $k=1$ and $v=w=0$, then we obtain the example from the introduction.
\end{enumerate}

\section{The Problem with higher order operators}
When dealing with higher order operators, the reasoning in the last section quickly gets out of hand. In the second order case, for example, we have the following problem:
As explained in the Appendix, the symbol of a second order differential pairings is a mapping from
$$\begin{array}{c}
\odot^{2}\g_{1}\otimes\mathbb{V}\otimes\mathbb{W}\\
\oplus\\
\g_{1}\otimes\mathbb{V}\otimes\g_{1}\otimes\mathbb{W}\\
\oplus\\
\mathbb{V}\otimes\odot^{2}\g_{1}\otimes\mathbb{W}
\end{array},$$
because all the term in here will have geometric weight $\omega_{1}+\omega_{2}+2$. Therefore there are 
$$2\times|\{\mathbb{E}\subset \odot^{2}\g_{1}\otimes\mathbb{V}\otimes\mathbb{W}\}|+|\{\mathbb{E}\subset \g_{1}\otimes\mathbb{V}\otimes\g_{1}\otimes\mathbb{W}\}|$$
unknowns corresponding to the terms which are second order in $V$, those which are second order in $W$ and those which are first order in both. However, there are
$$2\times|\{\mathbb{E}\subset \g_{1}\otimes\mathbb{V}\otimes\g_{1}\otimes\mathbb{W}\}|$$
obstruction terms, so it is not clear that we should obtain any pairings at all if there are more obstruction terms than unknowns. On $\mathbb{CP}_{n}$, for example, one
can look at all the pairings between
$$V=\;\xoo{w}{0}{0}\;...\;\oo{0}{0}=\mathcal{O}(w)\quad\text{and}\quad W=\;\xoo{1+v}{0}{0}\;...\;\oo{0}{1}$$
that land in
$$\xoo{v+w-2}{1}{0}\;...\;\oo{0}{0}.$$
The terms at our disposal are
$$f\nabla_{a}\nabla_{b} X^{b},\;(\nabla_{a}f)(\nabla_{b}X^{b}),\;(\nabla_{b}f)(\nabla_{a}X^{b}-\frac{1}{n}\nabla_{c}X^{c}\delta_{a}{}^{b}),\; X^{b}\nabla_{b}\nabla_{a}f$$
and there are four obstruction terms
$$f\Upsilon_{a}\nabla_{b}X^{b},\;f\Upsilon_{b}\nabla_{a}X^{b},\;(\nabla_{a}f)\Upsilon_{b}X^{b},\;(\nabla_{b}f)\Upsilon_{a}X^{b}.$$
So one might expect that only the zero paring would be invariant. But, somehow miraculously from this point of view, this is not the case and we obtain a one parameter family of invariant pairings spanned by
\begin{eqnarray*} 
X^{b}\nabla_{b}\nabla_{a}f&-&\frac{(w-1)(n+1)}{(v+n+1)n}\nabla_{a}f\nabla_{b}X^{b}\\
-\frac{w-1}{v+1}\nabla_{b}f(\nabla_{a}X^{b}-\frac{1}{n}\nabla_{c}X^{c}\delta_{a}{}^{b})&+&\frac{w(w-1)}{(v+1)(v+n+1)}f\nabla_{a}\nabla_{b} X^{b}.
\end{eqnarray*}

\section{Higher order pairings}
This section deals with $M$-th order bilinear invariant differential pairings. The strategy employed is to define a linear invariant differential mapping that includes an arbitrary irreducible homogeneous bundle in some other homogeneous bundle, called {\bf $M$-bundle} (which is in fact a {\bf tractor} {\bf bundle}, see~\cite{g}, p.~7), that encodes all the possible differential operators up to order $M$ emanating from this bundle. We will then tensor two of those $M$-bundles together and project onto irreducible components. First of all, we have to define the $M$-bundles:

\subsection{The $M$-module}
Let $\{\alpha_{i}\}_{i=0,...,n-1} $ be the simple roots of $\g$ with corresponding fundamental weights $\{\omega_{i}\}_{i=0,...,n-1}$, i.e.~$B(\omega_{i},\alpha_{j}^{\vee})=\delta_{i,j}$. One can order the simple roots in such a way that $\alpha_{0}$ is the distinguished simple root in $\g$ that makes a root space $\g_{\alpha}$ lie in $\g_{i}$ if and only if $i$ is the coefficient of $\alpha_{0}$ in the expression of $\alpha$ in simple roots.  We will define a representation $\mathbb{V}_{M}(\mathbb{E})$ of $\g$ that is induced from a finite dimensional irreducible representation $\mathbb{E}$ of $\g_{0}^{S}$ in the following way: $(\mathfrak{h}^{S})^{*}$ can be considered as a subspace of $\h^{*}$ in such a way that    
 $\{\alpha_{i}\}_{i=1,...,n-1} $ are the simple roots of $\g_{0}^{S}$ with corresponding fundamental weights $\{\omega_{i}\}_{i=1,...,n-1}$. The highest weight $\lambda$ of $\mathbb{E}^{*}$ can then be written as $\lambda=\sum_{i=1}^{n-1}a_{i}\omega_{i}$ with~$a_{i}\geq 0$. $\mathbb{V}_{M}(\mathbb{E})$ is defined to be the finite dimensional irreducible representation of $\g$ which is dual to the representation with highest weight $\Lambda=M\omega_{0}+\lambda\in\h$.
In the Dynkin diagram notation this is easily described. There is one node in the Dynkin diagram for $\g$ which denotes the simple root~$\alpha_{0}$. If we erase that node and adjacent edges, we obtain the Dynkin diagram for~$\g_{0}^{S}$. A finite dimensional irreducible representation $\mathbb{E}$ of $\g_{0}^{S}$ is denoted by writing non-negative integers associated to the highest weight of $\mathbb{E}^{*}$ over the nodes of this new diagram. $\mathbb{V}_{M}(\mathbb{E})$ is then denoted by writing those numbers over their corresponding nodes in the Dynkin diagram for $\g$ and in addition writing $M$ over the node that corresponds to~$\alpha_{0}$.

\subsubsection{Example}
The $|1|$-grading on $\g=\mathfrak{sl}_{n+1}\C$ as in Example~3.1.1. implies $\g_{0}^{S}\cong\mathfrak{sl}_{n}\C$, so for every representation 
$$\mathbb{E}=\;\oo{a_{1}}{a_{2}}\;...\;\oo{a_{n-2}}{a_{n-1}}$$
of $\g_{0}^{S}$ and every constant $M\geq 1$, we define
$$\mathbb{V}_{M}(\mathbb{E})=\;\ooo{M}{a_{1}}{a_{2}}\;...\;\oo{a_{n-2}}{a_{n-1}},$$
a representation of $\g$.

\begin{lemma}
As a $\p$ module $\mathbb{V}_{M}(\mathbb{E})$ has a composition series 
$$\mathbb{V}_{M}(\mathbb{E})=\mathbb{V}_{0}+\mathbb{V}_{1}+...+\mathbb{V}_{N},$$
so that $\g_{i}\mathbb{V}_{j}\subseteq \mathbb{V}_{i+j}$ and $\mathbb{V}_{0}\cong\mathbb{E}$ as a $\g_{0}^{S}$-module. Alternatively, one can look at this composition series as a splitting of $\g_{0}$ modules into eigenspaces for the action of the grading Element $E$. Thus $\mathbb{E}$ acquires the structure of a $\g_{0}$-module.
\end{lemma}
\begin{proof}
If $\Lambda$ is the highest weight of $\mathbb{V}_{M}(\mathbb{E})^{*}$, then $\mathbb{V}_{j}^{*}$ consists of those weight spaces, whose weight is of the form
$\Lambda-j\alpha_{0}-\sum_{i=1}^{n-1}k_{i}\alpha_{i}$, with $k_{i}\geq 0$. Therefore the action of $\g_{i}$ maps $\mathbb{V}^{*}_{i+j}$ to $\mathbb{V}^{*}_{j}$. Dually we obtain a mapping $\g_{i}: \mathbb{V}_{j}\rightarrow\mathbb{V}_{i+j}$. Note that $\mathbb{E}\cong\mathbb{V}_{0}=\mathbb{V}_{M}(\mathbb{E})/(\mathbb{V}_{1}+...+\mathbb{V}_{N})$ even acquires the structure of an irreducible $\p$-module that we can denote by a Dynkin diagram for $\p$: The integers over the uncrossed nodes correspond to the highest weight of $\mathbb{E}^{*}$ as a $\g_{0}^{S}$-module and $M$ is written over the crossed through node. 
\end{proof}


\begin{lemma}
There are $\g_{0}$ homomorphisms $\phi_{i}:\mathbb{V}_{i}\rightarrow\bigotimes^{i}\g_{1}\otimes\mathbb{E}$ that are injective for all $i$, have values in $\odot^{i}\g_{1}\otimes\mathbb{E}$ and define isomorphisms 
$$\mathbb{V}_{i}\cong\odot^{i}\g_{1}\otimes\mathbb{E}$$
for $0\leq i\leq M$.
\end{lemma}
\begin{proof}
The homomorphisms $\phi_{i}$ are constructed in ~\cite{bceg}, p.~655, with the help of Lie algebra cohomology and Kostant's version of the Bott-Borel-Weil theorem. The statement then follows from ~\cite{bceg}, p.~655, lemma~3.
\end{proof}

The next step is to look at a tensor product $\mathbb{V}_{M}(\mathbb{E})\otimes\mathbb{W}_{M}(\mathbb{F})$ and decompose it into irreducible $\g$-modules that themselves have composition series as $\p$-modules. The composition factors of all the irreducible components will then make up the composition factors of the tensor product. To be more precise, the composition series $\mathbb{V}_{M}(\mathbb{E})=\mathbb{V}_{0}+...+\mathbb{V}_{N_{1}}$ and $\mathbb{W}_{M}(\mathbb{F})=\mathbb{W}_{0}+...+\mathbb{W}_{N_{2}}$ induce a filtration on the tensor product:
$$\begin{array}{ccccccccc}
&&\;0-\text{th slot}\;&&\;1-\text{st slot}\;&&&&\;(N_{1}+N_{2})-\text{th slot}\\
\mathbb{V}_{M}(\mathbb{E})\otimes\mathbb{W}_{M}(\mathbb{F})&=&\mathbb{V}_{0}\otimes\mathbb{W}_{0}&+&
\begin{array}{c}
\mathbb{V}_{0}\otimes\mathbb{W}_{1}\\
\oplus\\
\mathbb{V}_{1}\otimes\mathbb{W}_{0}
\end{array}
&+&...&+&\mathbb{V}_{N_{1}}\otimes\mathbb{W}_{N_{2}}
\end{array}.$$
The $\g_{0}$ tensor product $\mathbb{V}_{i}\otimes\mathbb{W}_{j}$ can be decomposed into irreducible components by computing the irreducible components of the $\g_{0}^{S}$ tensor product $\mathbb{V}_{i}\otimes\mathbb{W}_{j}$. The geometric weight of all those components will be $\omega_{1}+\omega_{2}+i+j$, where $\omega_{1},\omega_{2}$ are the geometric weights of $\mathbb{E}$ and $\mathbb{F}$ respectively. Each of the irreducible components can be denoted by a Dynkin diagram for $\p$: The numbers over uncrossed nodes correspond to the highest weight of the irreducible component of the $\g_{0}^{S}$ tensor product $(\mathbb{V}_{i}\otimes\mathbb{W}_{j})^{*}$ that is to be denoted and the number over the crossed through node will make the geometric weight equal $\omega_{1}+\omega_{2}+i+j$. The next remark gives an estimate on those numbers which is of importance in the next proposition.

\subsubsection{Remark}
Let $\Lambda$ and $\lambda$ be defined as above and define $\mathbb{W}_{M}(\mathbb{F})$ analogously with highest weights $\Sigma$ and $\sigma$ of $\mathbb{W}_{M}(\mathbb{F})^{*}$ and $\mathbb{F}^{*}$ respectively.
All irreducible components in the $j$-th slot of $\mathbb{V}_{M}(\mathbb{E})\otimes\mathbb{W}_{M}(\mathbb{F})$ are dual to representations of highest weights of the form $\mu=\Lambda+\Sigma-j\alpha_{0}-\sum_{i=1}^{n-1}k_{i}\alpha_{i}$, so the number over the crossed through node will be
$$(\mu,\alpha_{0}^{\vee})=2M-2j-\sum_{i=1}^{n-1}k_{i}(\alpha_{i},\alpha_{0}^{\vee})\geq 2(M-j),$$
since $(\alpha_{i},\alpha_{0}^{\vee})\leq 0$ for all $i=1,...,n-1$.

\begin{proposition}
Let $\mathbb{E}$ and $\mathbb{F}$ be two finite dimensional irreducible representations of  $\g_{0}^{S}$.
If $l\leq M$, then for every irreducible component $\mathbb{H}$ of the $\g_{0}^{S}$ tensor product $\odot^{l}\g_{1}\otimes\mathbb{E}\otimes\mathbb{F}$, there is a $\p$-module projection
$$\mathbb{V}_{M}(\mathbb{E})\otimes\mathbb{W}_{M}(\mathbb{F})\rightarrow \mathbb{H},$$
where $\mathbb{H}$ has acquired the structure of an irreducible $\p$-module with geometric weight $\omega_{1}+\omega_{2}+l$.
\end{proposition}\begin{proof}
As $\p$-modules, the $M$-modules associated to $\mathbb{E}$ and $\mathbb{F}$ have composition series
$$\mathbb{V}_{M}(\mathbb{E})=\mathbb{E}+\g_{1}\otimes \mathbb{E}+\odot^{2}\g_{1}\otimes \mathbb{E}+...+\odot^{M}\g_{1}\otimes \mathbb{E}+\mathbb{V}_{M+1}+...+\mathbb{V}_{N_{1}}$$
and
$$\mathbb{W}_{M}(\mathbb{F})=\mathbb{F}+\g_{1}\otimes \mathbb{F}+\odot^{2}\g_{1}\otimes \mathbb{F}+...+\odot^{M}\g_{1}\otimes \mathbb{F}+\mathbb{W}_{M+1}+...+\mathbb{W}_{N_{2}}\;.$$
Therefore the tensor product $\mathbb{V}_{M}(\mathbb{E})\otimes\mathbb{W}_{M}(\mathbb{F})$ has a composition series
$$\mathbb{V}_{M}(\mathbb{E})\otimes\mathbb{W}_{M}(\mathbb{F})=\mathbb{E}\otimes\mathbb{F}+\begin{array}{c}
\mathbb{E}\otimes\g_{1}\otimes \mathbb{F}\\
\oplus\\
\g_{1}\otimes \mathbb{E}\otimes \mathbb{F}
\end{array}
+\begin{array}{c}
\mathbb{E}\otimes\odot^{2}\g_{1}\otimes \mathbb{F}\\
\oplus\\
\g_{1}\otimes \mathbb{E}\otimes\g_{1}\otimes \mathbb{F}\\
\oplus\\
\odot^{2}\g_{1}\otimes \mathbb{E}\otimes\mathbb{F}
\end{array}+...\;.$$
Every irreducible component $\mathbb{G}$ of $\mathbb{E}\otimes\mathbb{F}$ (as $\g_{0}^{S}$-modules) corresponds to an irreducible component $\mathbb{U}$ in $\mathbb{V}_{M}(\mathbb{E})\otimes\mathbb{W}_{M}(\mathbb{F})$ (as $\g$-modules) that has a composition series that starts with $\mathbb{G}$ and then continues with $\g_{1}\otimes\mathbb{G}+\odot^{2}\g_{1}\otimes \mathbb{G}+...\;$. We will say that the composition series is {\bf predictable} up to the $x$-th slot, if $\mathbb{U}_{j}\cong \odot^{j}\g_{1}\otimes\mathbb{G}$ for all $j\leq x$, as $\g_{0}$-modules. Using Lemma~4, we know that the composition series  of $\mathbb{U}$ is predictable up to the $x$-th slot if the number over the crossed through node in $\mathbb{G}$ is~$x$.
\par
Removing all those composition series corresponding to irreducible components of $\mathbb{E}\otimes\mathbb{F}$ from the composition series of $\mathbb{V}_{M}(\mathbb{E})\otimes\mathbb{W}_{M}(\mathbb{F})$ leaves nothing in the zeroth slot, exactly one copy of $\mathbb{E}\otimes\g_{1}\otimes \mathbb{F}$ in the first slot, one copy of each 
$\g_{1}\otimes \mathbb{E}\otimes\g_{1}\otimes \mathbb{F}$ and $ \mathbb{E}\otimes\odot^{2}\g_{1}\otimes\mathbb{F}$
in the second slot and so forth. Therefore the next irreducible components of $\mathbb{V}_{M}(\mathbb{E})\otimes\mathbb{W}_{M}(\mathbb{F})$ all have a composition series that starts with an irreducible component of $\mathbb{E}\otimes\g_{1}\otimes \mathbb{F}$. Removing those again leaves nothing in the first two slots, exactly one copy of
$\mathbb{E}\otimes\odot^{2}\g_{1}\otimes \mathbb{F}$ in the second slot and so forth.
Hence the next irreducible components of $\mathbb{V}_{M}(\mathbb{E})\otimes\mathbb{W}_{M}(\mathbb{F})$ correspond to irreducible components of $\mathbb{E}\otimes\odot^{2}\g_{1}\otimes \mathbb{F}$. This argument is correct as long as all the compositions series are predictable. This is the case exactly up to the $M$-th slot:
\par
In the $l\leq M$-th slot of $\mathbb{V}_{M}(\mathbb{E})\otimes\mathbb{W}_{M}(\mathbb{F})$ the lowest number over a cross is bigger or equal to
$2(M-l)$ (Remark~6.1.2.). Some of the factors in here correspond to irreducible components of $\mathbb{V}_{M}(\mathbb{E})\otimes\mathbb{W}_{M}(\mathbb{F})$  as
$\g$-representations that themselves have a composition series that is predictable up the the $2(M-l)$-th slot, which corresponds in the big composition series to the $(2M-l)$-th slot. So the argument above is correct for $l\leq M$. There could be (and in general this happens) more irreducible components of $\mathbb{V}_{M}(\mathbb{E})\otimes\mathbb{W}_{M}(\mathbb{F})$, but those correspond to higher order pairings.  The mapping is defined by first projecting onto the correct irreducible component of $\mathbb{V}_{M}(\mathbb{E})\otimes\mathbb{W}_{M}(\mathbb{F})$ and then projecting onto the first composition factor in the composition series, which will, as a $\g_{0}^{S}$-module, be isomorphic to $\mathbb{H}$. The $\p$-module structure is derived as in Lemma~3. 
\end{proof}

\subsubsection{Remark}
For every $k\in \Z$, let $\mathcal{O}(k-M)$ be the one dimensional $\p$-module which is dual to the representation of highest weight $(k-M)\omega_{0}$. In the Dynkin diagram notation this corresponds to having $k-M$ over the crossed through node and zeros elsewhere. We will write
$\mathbb{V}(k-M)$ for the tensor product $\mathbb{V}\otimes\mathcal{O}(k-M)$ for every $\p$-module $\mathbb{V}$ and remark that this procedure only changes the geometric weight of $\mathbb{V}$. Hence if $\mathbb{V}$ has a composition series, then the composition series of $\mathbb{V}(k-M)$ is obtained by tensoring each factor with $\mathcal{O}(k-M)$. Furthermore $\mathbb{V}_{0}(k-M)$ is dual to a representation of highest weight $k\omega_{0}+\sum_{i=1}^{n-1}a_{i}\omega_{i}$ and by choosing $\mathbb{E}$ and $k$ correctly, we can write every finite dimensional irreducible $\p$-module as $\mathbb{V}_{0}(k-M)$ of some module $\mathbb{V}_{M}(\mathbb{E})(k-M)$.

\par
\vspace{0.2cm}
The idea is now to define invariant linear differential mappings 
$$V_{0}(k-M)\rightarrow V_{M}(\mathbb{E})(k-M)\quad\text{and}\quad W_{0}(l-M)\rightarrow W_{M}(\mathbb{F})(l-M).$$
We can then tensor 
$$\mathbb{V}_{M}(\mathbb{E})(k-M)\otimes \mathbb{W}_{M}(\mathbb{E})(l-M) =(\mathbb{V}_{M}(\mathbb{E})\otimes \mathbb{W}_{M}(\mathbb{F}))(k+l-2M)$$
together and project onto the first composition factor of every irreducible component of $(\mathbb{V}_{M}(\mathbb{E})\otimes \mathbb{W}_{M}(\mathbb{F}))(k+l-2M)$ just as in Proposition~1, with the only difference that the geometric weight of $\mathbb{H}$ will be different. This is clearly a bilinear invariant differential pairing between sections of $V_{0}(k-M)$ and $W_{0}(l-M)$. In order to do this in the flat homogeneous case $G/P$ we can use the following three theorems:

\begin{theorem}
Invariant linear differential operators between sections of homogeneous bundles over a flag manifold $G/P$ are in one-to-one correspondence with $\g$-module homomorphisms of induced modules.
\end{theorem}\begin{proof}
This theorem is proved in a straightforward manner in~\cite{er}, p.~212. It may be noted that the theorem is usually stated in terms of generalized Verma modules (see~~\cite{be}, p.~164), the statement, however, remains true for induced modules with identical proof.
\end{proof}

\begin{theorem}
If $M_{\p}(\mathbb{V}_{0}(k-M))$ has distinct central character from the generalized Verma modules associated to all the other composition factors of $\mathbb{V}_{M}(\mathbb{E})(k-M)$, then it can be canonically split off as a direct summand of $M_{\p}(\mathbb{V}_{M}(\mathbb{E})(k-M))$.
\end{theorem}
\begin{proof}
The composition series $\mathbb{V}_{M}(\mathbb{E})(k-M)=\mathbb{V}_{0}(k-M)+...+\mathbb{V}_{N}(k-M)$ induces a composition series
$$(\mathbb{V}_{M}(\mathbb{E})(k-M))^{*}=(\mathbb{V}_{N}(k-M))^{*}+...+(\mathbb{V}_{0}(k-M))^{*}$$
of the dual representation. Since the functor that associates to every $\p$-module $\mathbb{V}^{*}$ the corresponding induced module $\mathfrak{U}(\g)\otimes_{\mathfrak{U}(\p)}\mathbb{V}^{*}$ is exact (see ~\cite{v}, p.~303, lemma~6.1.6), we have a filtration 
$$M_{\p}(\mathbb{V}_{M}(\mathbb{E})(k-M))=M_{\p}(\mathbb{V}_{N}(k-M))+...+M_{\p}(\mathbb{V}_{0}(k-M))$$
that induces an injection $M_{\p}(\mathbb{V}_{0}(k-M))\hookrightarrow M_{\p}(\mathbb{V}_{M}(\mathbb{E})(k-M))$.
The weight spaces of $M_{\p}(\mathbb{V}_{M}(\mathbb{E})(k-M))$ can be grouped in terms of central character, so the projection $M_{\p}(\mathbb{V}_{M}(\mathbb{E})(k-M))\rightarrow M_{\p}(\mathbb{V}_{0}(k-M))$ may be defined by projecting onto the joint eigenspace of the central character of $M_{\p}(\mathbb{V}_{0}(k-M))$. Since central character is preserved under the action of $\g$, this projection is indeed a $\g$-module homomorphism and provides a $\g$-module splitting  of $M_{\p}(\mathbb{V}_{M}(\mathbb{E})(k-M))$.
\end{proof}

\begin{theorem}[Harish-Chandra]
Two generalized Verma modules have the same central character if and only if their highest weights are related by the affine action of the Weyl group of $\g$.
\end{theorem}
\begin{proof}
A proof of this theorem can, for example, be found in ~\cite{h}, p.~130, theorem~23.3.
\end{proof}
\par
These three theorems combined are the backbone of the {\bf Jantzen}-{\bf Zuckermann} 
{\bf translation} {\bf functor} as used in ~\cite{es} and ~\cite{er}.  In principle they can be used to define invariant bilinear differential pairings for every homogeneous space $G/P$ with an AHS structure. One only has to exclude weights, i.e.~values of $k$, for which the central character of $M_{\p}(\mathbb{V}_{0}(k-M))$ is the same as the central character of a generalized Verma module associated to another
composition factor of $\mathbb{V}_{M}(\mathbb{E})(k-M)$. A trivial case is $k=M$, because all the weight spaces of $\mathbb{V}_{M}(\mathbb{E})$ apart from the highest weight space, which lies in $\mathbb{V}_{0}$, have weights $\mu$ so that
$$\Vert \Lambda+\rho_{\g}\Vert^{2} >\Vert \mu+\rho_{\g}\Vert^{2},$$
with $\rho_{\g}=\frac{1}{2}\sum_{\alpha\in\Delta^{+}(\g)}\alpha$ (see~~\cite{h}, p.~114, proposition~21.4 and p.~71, lemma~13.4).
Since the Weyl group acts by isometries, this implies that $M_{\p}(\mathbb{V}_{0})$ has distinct central character from the generalized Verma modules associated to all the other composition factors of $\mathbb{V}_{M}(\mathbb{E})$. The pairings that we obtain via our construction are then the flat analogues of the parings $\sqcup_{\eta}$ as defined in
~\cite{cd}, p.~13, theorem~3.6.

\section{Higher order pairings for $\mathbb{CP}_{n}$}
In this section we will work exclusively on complex projective space $\mathbb{CP}_{n}$ with $\g=\mathfrak{sl}_{n+1}\C$ and the conventions described above. A Dynkin diagram therefore stands for four things: An irreducible representation of $\p$, the corresponding irreducible homogeneous vector bundle, its sections and the generalized Verma module associated to the representation. In every case it should be clear which meaning we refer to and sometimes it is convenient that two meanings are denoted at the same time.
\par
\begin{definition}
\em For every representation 
$$\mathbb{E}=\;\oo{a_{1}}{a_{2}}\;...\;\oo{a_{n-2}}{a_{n-1}}$$
of $\g_{0}^{S}=\mathfrak{sl}_{n}\C$ and every constant $M\geq 1$ we define
$$\mathbb{V}_{M}(\mathbb{E})=\;\ooo{M}{a_{1}}{a_{2}}\;...\;\oo{a_{n-2}}{a_{n-1}},$$
a representation of $\g$, which we also denote by
$$\mathbb{V}_{M}(\mathbb{E})=(0,b_{0},b_{1},b_{2},....,b_{n-1})=\left(0,M,a_{1}+M,a_{1}+a_{2}+M,...,\sum_{i=1}^{n-1}a_{i}+M\right).$$
When referring to a representation of $\p$, we will use the notation $(a|b,c,...,d,e,f)$ for $\;\xo{b-a}{c-b}\;...\;\oo{e-d}{f-e}.$
This is important whenever we want to describe the action of the Weyl group $\mathcal{W}$ on the weight, because $\mathcal{W}\cong\mathbb{S}_{n+1}$ and it acts  by permutation. 
\end{definition}
\par
\vspace{0.2cm} 
The $\g$-module $\mathbb{V}_{M}(\mathbb{E})$ has, as a $\p$-module, a composition series
$$\mathbb{V}_{M}(\mathbb{E})=\mathbb{V}_{0}+\mathbb{V}_{1}+\mathbb{V}_{2}+...+\mathbb{V}_{N},$$
where each $\mathbb{V}_{i}$ decomposes into a direct sum of irreducible $\p$-modules and 
$$\mathbb{V}_{0}=\xooo{M}{a_{1}}{a_{2}}{a_{3}}\;...\;\oo{a_{n-2}}{a_{n-1}}.$$
We may tensor this composition series by $\mathcal{O}(k-M)$ to obtain
$$\mathbb{V}_{0}(k-M)= \xooo{k}{a_{1}}{a_{2}}{a_{3}}\;...\;\oo{a_{n-2}}{a_{n-1}}.$$
This is the $\p$-module that we are interested in and we want to define a mapping 
$$V_{0}(k-M)\rightarrow V_{M}(\mathbb{E})(k-M)$$
using the theorems from the last section. Hence we have to make sure that the generalized Verma modules associated to all the irreducible composition factors of $\mathbb{V}_{M}(\mathbb{E})(k-M)$ have a central character which is different from the central character of~$M_{\p}(\mathbb{V}_{0}(k-M))$.

\subsubsection{Remark}
In the case of  $\mathbb{CP}_{n}$, lemma~4 can be proved directly using Pierie's formula, as in ~\cite{fh},~p.~225, for the tensor product $\odot^{l}\g_{1}\otimes\mathbb{E}$ and the branching rules for restrictions of representations of $\mathfrak{sl}_{n+1}\C$ to $\mathfrak{sl}_{n}\C$ as in ~\cite{gw},~p.~350. The upshot of this procedure is that 
$\mathbb{V}_{l}(k-M)$ consists of terms $(M-k+l|\tilde{b}_{0},\tilde{b}_{1},...,\tilde{b}_{n-1})$ that {\bf interlace} $(M-k|b_{0},b_{1},....,b_{n-1})$, i.e.
$$0\leq\tilde{b}_{0}\leq b_{0}\leq \tilde{b}_{1}\leq b_{1}\leq\tilde{b}_{2}\leq b_{2}\leq...\leq\tilde{b}_{n-1}\leq b_{n-1}$$
and  $\sum_{i=0}^{n-1}b_{i}-\sum_{i=0}^{n-1}\tilde{b}_{i}=l$.
We can also see that $N=\sum_{i=1}^{n-1}a_{i}+M$, because for $l>N$ it is not possible for any $(M-k+l|\tilde{b}_{0},\tilde{b}_{1},...,\tilde{b}_{n-1})$ to interlace~$(M-k|b_{0},b_{1},...,b_{n-1})$.

\begin{proposition}
The only irreducible components of $\mathbb{V}_{l}(k-M)$ that can induce generalized Verma modules with the same central character as $M_{\p}(\mathbb{V}_{0}(k-M))$
are the ones that are of the form
$$(M-k+l|b_{0},b_{1},...,b_{j-1},b_{j}-l,b_{j+1},...,b_{n-1}),$$
for $j=0,1,...,n-1$. If $j\in\{1,...,n-1\}$, then this is only allowed for $a_{j}\geq l$ and if $j=0$, then this is only allowed for $l\leq M$.
In that case the generalized Verma module has the same central character as $M_{\p}(\mathbb{V}_{0}(k-M))$ if and only if
$$k=-\left(\sum_{i=1}^{j}a_{i}+j-l+1\right).$$
For $j=0$, this condition reads $k=l-1$.
\end{proposition}\begin{proof}
Using remark~7.0.4, we know that an arbitrary irreducible component $\mathbb{V}_{l,v}(k-M)$ of $\mathbb{V}_{l}(k-M)$ has to be of the form $(M-k+l|\tilde{b}_{0},...,\tilde{b}_{n-1})$, so that $(M-k+l|\tilde{b}_{0},...,\tilde{b}_{n-1})$ interlaces~$(M-k|b_{0},...,b_{n-1})$.  Let us assume that there are at least two integers $0\leq i<j\leq n-1$, such that $\tilde{b}_{i}<b_{i}$ and~$\tilde{b}_{j}<b_{j}$. We can assume that $i$ is the smallest integer with this property and that $j$ is the biggest integer with this property.
\par
Theorem~4 implies that the central characters of $M_{\p}(\mathbb{V}_{0}(k-M))$ and $M_{\p}(\mathbb{V}_{l,v}(k-M))$ are identical if and only if there is an element in the Weyl group, i.e.~a permutation, that
maps $(M-k+l|\tilde{b}_{0},...,\tilde{b}_{n-1})+\rho_{\g}$ to $(M-k|b_{0},...,b_{n-1})+\rho_{\g}$. Using $\rho_{\g}=(1,2,...,n,n+1)$, we obtain the condition that the two sets
$$\{M-k+1,b_{0}+2,b_{1}+3,...,b_{i}+i+2,...,b_{j}+j+2,...,b_{n-1}+n+1\}$$
and
$$\{M-k+l+1,\tilde{b}_{0}+2,\tilde{b}_{1}+3,...,\tilde{b}_{i}+i+2,...,\tilde{b}_{j}+j+2,...,\tilde{b}_{n-1}+n+1\}$$
have to be equal.
This is equivalent to
$$\{M-k+1,b_{i}+i+2,...,b_{j}+j+2\}=\{M-k+l+1,\tilde{b}_{i}+i+2,...,\tilde{b}_{j}+j+2\},$$
where the sets contain all those $b_{m}+m+2$, resp. $\tilde{b}_{m}+m+2$, for which~$\tilde{b}_{m}\not=b_{m}$. Furthermore, leaving out $M-k+1$, all numbers in the first set are increasing from left to right.
Since $\tilde{b}_{i}<b_{i}$, $\tilde{b}_{i}+i+2$ is smaller than the second entry in the first set and therefore smaller than everything but the first entry, i.e.~we must
have~$\tilde{b}_{i}+i+2=M-k+1$. 
Moreover $\tilde{b}_{j}<b_{j}$ implies that there has to be an integer $m<j$, so that  
\begin{eqnarray*}
&\tilde{b}_{j}+j+2=b_{m}+m+2
\Rightarrow&\tilde{b}_{j}+j= b_{m}+m.
\end{eqnarray*}
This is not possible, because $\tilde{b}_{j}\geq b_{m}$ and~$j>m$. That proves the first claim.
\par
Let us now assume that $\mathbb{V}_{l,v}(k-M)=(k-M+l|b_{0},b_{1},...,b_{j-1},b_{j}-l,b_{j+1},...,b_{n-1})$. In this case $M_{\p}(\mathbb{V}_{l,v}(k-M))$ has the same central character as $M_{\p}(\mathbb{V}_{0}(k-M))$ if and only if
$$\{M-k+l+1,b_{j}-l+j+2\}\;=\;\{M-k+1,b_{j}+j+2\},$$
which is equivalent to $k=-b_{j}+M-j+l-1=-\left(\sum_{i=1}^{j}a_{i}+j-l+1\right)$.
\end{proof}

\begin{definition}
\rm Let $\mathfrak{h}^{S}$ denote the Cartan subalgebra of~$\g_{0}^{S}$. Then we have
$$(\mathfrak{h}^{S})^{*}=\C\langle L_{1},...,L_{n}\rangle/(L_{1}+...+L_{n}=0),$$
where $L_{i}(H_{j})=\delta_{i,j}$ and $H_{j}$ denotes the matrix which has a one in the $j$-th diagonal entry and zeros elsewhere as an element in~$\mathfrak{sl}_{n}\C$.
\end{definition}

\begin{proposition}
If
$$k=-\left(\sum_{i=1}^{j}a_{i}+j-l+1\right),$$
then there is a $l$-th order invariant linear differential operator  
$$\xooo{k}{a_{1}}{a_{2}}{a_{3}}\;...\;\oo{a_{n-2}}{a_{n-1}}\;\rightarrow \xoo{k-l}{a_{1}}{a_{2}}\;...\oo{a_{j}-l}{a_{j+1}+l}\;...\;\oo{a_{n-2}}{a_{n-1}}.$$
\end{proposition}
\begin{proof}
As proved in ~\cite{css},~p.~65, corollary~5.3, the condition for this operator to be invariant is
$$\omega=(\alpha+L_{n-j},\rho)-\frac{1}{2}(l-1)(|\alpha|^{2}+1)-(-L_{n-j},\tilde{\lambda}),$$
where $\omega=-\frac{1}{n+1}\left(nk+\sum_{i=1}^{n-1}(n-i)a_{i}\right)$ is the geometric weight of $\xo{k}{a_{1}}\;...\;\oo{a_{n-2}}{a_{n-1}}$, $(.,.)$ is the normalized Killing form as in the previous sections and $\alpha=-L_{n}$ is the highest weight of~$\g_{1}$. Moreover $\rho=\rho_{\mathfrak{sl}_{n}\C}=\sum_{i=1}^{n-1}(n-i)L_{i}$, $|\alpha|^{2}=(\alpha,\alpha)$ and $\tilde{\lambda}=\sum_{i=1}^{n-1}\lambda_{i}L_{i}$ (we can always assume that $\lambda_{n}=0$, which implies $\lambda_{n-j}=\sum_{i=1}^{j}a_{j}$) is the highest weight 
of~$\mathbb{E}$. 
Using
\begin{eqnarray*}
(\alpha+L_{n-j},\rho)&=&\frac{nj}{n+1},\\
|\alpha|^{2}&=&\frac{n-1}{n+1},\\
(L_{n-j},\tilde{\lambda})&=&\frac{n\lambda_{n-j}-\sum_{i=1}^{n}\lambda_{i}}{n+1}
\end{eqnarray*}
and the formula for $\omega$ from above, we see that
$$\omega=(\alpha+L_{n-j},\rho)-\frac{1}{2}(l-1)(|\alpha|^{2}+1)-(-L_{n-j},\tilde{\lambda})\Leftrightarrow k=-\left(\sum_{i=1}^{j}a_{i}+j-l+1\right).$$
Note that these calculations for $j\in\{1,...,n-1\}$ make only sense if $a_{j}\geq l$. If $j=0$, then $l$ may be arbitrary.
\end{proof}

The problem is, when we look at $M$-th order pairings, we do not really want to exclude weights that correspond to operators that have a higher order. The following 
lemma excludes such a situation at the cost of a restriction on the integers~$a_{i}$.

\begin{lemma}
Let $M\geq \mathrm{max}_{i}\{a_{i}\}$, then no weights have to be excluded for~$l> M$:
\end{lemma}
\begin{proof}
As discussed earlier, an irreducible component in $\mathbb{V}_{l}(k-M)$ that induces a generalized Verma module with the same central character as $M_{\p}(\mathbb{V}_{0}(k-M))$ can only arise by taking
$$(M-k|b_{0},b_{1},...,b_{n-1})$$
and subtracting $l$ from one of the $b_{i}$ to obtain
$$(M-k+l|\tilde{b}_{0},\tilde{b}_{1},...,\tilde{b}_{n-1}),$$
so that $(M-k+l|\tilde{b}_{0},...,\tilde{b}_{n-1})$ interlaces~$(M-k|b_{0},...,b_{n-1})$.
But $b_{i}-b_{i-1}=a_{i}\leq M< l\;\forall\;i=1,...,n-1$, so subtracting $l$ from any $b_{i}$, $i\geq 1$, leads to $\tilde{b}_{i}=b_{i}-l<b_{i-1}$, which is not allowed. 
Subtracting $l$ from $b_{0}$ leaves $\tilde{b}_{0}=M-l<0$, which is also not allowed.
Therefore all terms in $\mathbb{V}_{l}(k-M)$, for $l>M$, induce a generalized Verma module that has a central character which is different from the one of~$M_{\p}(\mathbb{V}_{0}(k-M))$.
\end{proof}

\subsubsection{Examples}
\begin{enumerate}
\item
Let us look at symmetric  two tensors of projective weight $v$, i.e.~sections of $\odot^{2}TM\otimes\mathcal{O}(v)$ for~$M=2$:
\begin{eqnarray*}
\oo{2}{0}\;...\;\oo{0}{2}(v)\;&=&\;\xo{2+v}{0}\;...\;\oo{0}{2}\;+\begin{array}{c}
\;\xo{1+v}{0}\;...\;\oo{0}{1}\;\\
\oplus\\
\;\xoo{v}{1}{0}\;...\;\oo{0}{2}\;
\end{array}
+\begin{array}{c}
\;\xo{v}{0}\;...\;\oo{0}{0}\;\\
\oplus\\
\;\xoo{v-1}{1}{0}\;...\;\oo{0}{1}\;\\
\oplus\\
\;\xoo{v-2}{2}{0}\;...\;\oo{0}{2}\;
\end{array}\\
&&+\begin{array}{c}
\;\xoo{v-2}{1}{0}\;...\;\oo{0}{0}\;\\
\oplus\\
\;\xoo{v-3}{2}{0}\;...\;\oo{0}{1}\;
\end{array}+
\;\xoo{v-4}{2}{0}\;...\;\oo{0}{0}\;.\end{eqnarray*}
The weights to exclude are
\begin{enumerate}
\item
$v=-2,-(n+3)$ which correspond to invariant first order operators $\nabla_{a}V^{bc}-\frac{2}{n+1}\delta_{a}{}^{(b}\nabla_{d}V^{c)d}$ and $\nabla_{a}V^{ab}$ respectively;
\item
$v=-1,-(n+2)$ which correspond to invariant second order operators $\nabla_{a}\nabla_{b}V^{cd}-\mathrm{trace}$ and $\nabla_{a}\nabla_{b}V^{ab}$ respectively.
\end{enumerate}
\item
Another example for vector fields of projective weight $v$, i.e. sections of $TM\otimes\mathcal{O}(v)$, with~$M=1$:
$$\;\oo{1}{0}\;...\;\oo{0}{1}(v)\;=\;\xo{1+v}{0}\;...\;\oo{0}{1}\;+\begin{array}{c}
\;\xo{v}{0}\;...\;\oo{0}{0}\;\\
\oplus\\
\;\xoo{v-1}{1}{0}\;...\;\oo{0}{1}\;
\end{array}
+\;\xoo{v-2}{1}{0}\;...\;\oo{0}{0}\;.$$
The weights to exclude are $v=-1,-(n+1)$ corresponding to invariant first order operators $\nabla_{a}V^{b}-\frac{1}{n}\delta_{a}{}^{b}\nabla_{c}V^{c}$ and $\nabla_{a}V^{a}$
respectively.
\item
The last example deals with weighted functions and a general~$M$:
\begin{eqnarray*}\;\oo{M}{0}\;...\;\oo{0}{0}(w-M)\;&=&\;\xo{w}{0}\;...\;\oo{0}{0}\;+\;\xoo{w-2}{1}{0}\;...\;\oo{0}{0}\;+\;\xoo{w-4}{2}{0}\;...\;\oo{0}{0}\;\\
&&+...+\;\xoo{w-2M}{M}{0}\;...\;\oo{0}{0}\;.\end{eqnarray*}
The weights to exclude are $w=0,1,...M-1$ corresponding to invariant operators $\underbrace{\nabla_{a}...\nabla_{c}}_{w+1}f$ respectively.
\end{enumerate}

To state the main theorem, we have to define precisely what we mean by excluded weights.

\begin{definition}
\rm Let $\;\xoo{k}{a_{1}}{a_{2}}\;...\;\oo{a_{n-2}}{a_{n-1}}\;$ be a representation of $\p$. Then the {\bf excluded} {\bf weights} up to order $M$ consist of 
all $k$ such that there is a $1\leq l\leq M$ and a $0\leq j\leq n-1$ with
$$k=-\left(\sum_{i=1}^{j}a_{i}+j-l+1\right)\quad\text{and}\; a_{j}\geq l.$$
For $j=0$, the excluded weights are $k=l-1$ for~$1\leq l\leq M$.
\end{definition}

\begin{theorem}[Main Result 2]
Let $\;\xoo{k}{a_{1}}{a_{2}}\;...\;\oo{a_{n-2}}{a_{n-1}}\;$ and $\;\xoo{m}{b_{1}}{b_{2}}\;...\;\oo{b_{n-2}}{b_{n-1}}\;$ be irreducible homogeneous bundles on $\mathbb{CP}_{n}$.
If $M\geq \max_{i}\{a_{i},b_{i}\}$ and $k$ and $m$ are not equal to one of the excluded weights up to order $M$, then there exists an $r$ parameter family of $M$-th order bilinear invariant differential pairings
$$ \;\xoo{k}{a_{1}}{a_{2}}\;...\;\oo{a_{n-2}}{a_{n-1}}\;\times\;\xoo{m}{b_{1}}{b_{2}}\;...\;\oo{b_{n-2}}{b_{n-1}}\;\rightarrow \;\xoo{s}{c_{1}}{c_{2}}\;...\;\oo{c_{n-2}}{c_{n-1}},\;$$
where $r$ is the multiplicity of $\;\oo{c_{1}}{c_{2}}\;...\;\oo{c_{n-2}}{c_{n-1}}\;$ in 
$$\odot^{M}\g_{1}\otimes\;\oo{a_{1}}{a_{2}}\;...\;\oo{a_{n-2}}{a_{n-1}}\;\otimes\;\oo{b_{1}}{b_{2}}\;...\;\oo{b_{n-2}}{b_{n-1}}\;.$$
Excluded weights correspond to invariant linear differential operators of order $\leq M$ emanating from the bundles in question.
\end{theorem}
\begin{proof}
If $M\geq \max_{i}\{a_{i},b_{i}\}$ and $k$ and $m$ are not equal to one of the excluded weights up to order $M$, we can use lemma 5, proposition 2, theorem 3 and theorem 2 to define
invariant differential operators that take $\;\xoo{k}{a_{1}}{a_{2}}\;...\;\oo{a_{n-2}}{a_{n-1}}\;$ and $\;\xoo{m}{b_{1}}{b_{2}}\;...\;\oo{b_{n-2}}{b_{n-1}}\;$ into their $M$-bundles.
Then we decompose the tensor product of the $M$-bundles as described in proposition 1 and project onto the first composition factor of each of the irreducible components.
That also yields all the invariant pairings of order smaller than $M$, but we may have to exclude more weights than necessary. Moreover there cannot be more invariant pairings, because then one would be able to find a linear combination of all those pairings that does not involve the highest order terms ($M$ derivatives) in sections of one of the bundles. But obstruction terms involving $M-1$ derivatives in the sections of that bundle and one $\Upsilon$-term would therefore only occur in $\odot^{M-1}\g_{1}\otimes \mathbb{E}\otimes\g_{1}\otimes \mathbb{F}$ (if $\mathbb{E}$ and $\mathbb{F}$ denote the corresponding $\g_{0}^{S}$-modules as before) and one would not be able to eliminate them, because no operator in the formula is invariant. The last statement follows from Proposition~3.
\end{proof}

\subsubsection{Example}
Let us carry out the described construction for first order pairings between weighted $2$-forms and weighted vector fields on~$\mathbb{CP}_{4}$. The corresponding $M$ bundles have composition series
$$\;\oooo{1}{0}{1}{0}\;=\;\xooo{1}{0}{1}{0}\;+\begin{array}{c}
\;\xooo{0}{0}{0}{1}\;\\
\oplus\\
\;\xooo{-1}{1}{1}{0}\;
\end{array}
+\;\xooo{-2}{1}{0}{1}\;$$
and
$$\;\oooo{1}{0}{0}{1}\;=\;\xooo{1}{0}{0}{1}\;+\begin{array}{c}
\;\xooo{0}{0}{0}{0}\;\\
\oplus\\
\;\xooo{-1}{1}{0}{1}\;
\end{array}
+\;\xooo{-2}{1}{0}{0}\;.$$
If we tensor these together, we obtain a composition series
$$
\;\oooo{1}{0}{1}{0}\;\otimes\;\oooo{1}{0}{0}{1}\; =
\left(\begin{array}{c}
\;\xooo{2}{0}{1}{1}\;\\
\oplus\\
\;\xooo{2}{1}{0}{0}\;
\end{array}\right)
+\left(\begin{array}{ccc}
4\times\;\xooo{1}{0}{1}{0}\;&\oplus&2\times\;\xooo{0}{1}{1}{1}\;\\
&\oplus&\\
2\times\;\xooo{0}{2}{0}{0}\;&\oplus&2\times\;\xooo{1}{0}{0}{2}\;
\end{array}\right)
$$
$$
+\left(\begin{array}{ccc}
3\times\;\xooo{-1}{1}{0}{2}\;&\oplus&6\times\;\xooo{-1}{1}{1}{0}\;\\
&\oplus&\\
5\times\;\xooo{0}{0}{0}{1}\;&\oplus&\;\xooo{-2}{2}{1}{1}\;\\
&\oplus&\\
\;\xooo{-1}{0}{2}{1}\;&\oplus&\;\xooo{-2}{3}{0}{0}\;
\end{array}\right)
+\left(\begin{array}{ccc}
5\times\;\xooo{-2}{1}{0}{1}\;&\oplus&2\times\;\xooo{-1}{0}{0}{0}\;\\
&\oplus&\\
2\times\;\xooo{-3}{2}{1}{0}\;&\oplus&2\times\;\xooo{-2}{0}{2}{0}\;\\
&\oplus&\\
\;\xooo{-3}{2}{0}{2}\;&\oplus&\;\xooo{-2}{0}{1}{2}\;
\end{array}\right)
$$
$$
+\left(\begin{array}{c}
\;\xooo{-4}{2}{0}{1}\\
\oplus\\
\;\xooo{-3}{0}{1}{1}\;\\
\oplus\\
\;\xooo{0}{1}{0}{0}\;
\end{array}\right).
$$
This composition series can be split up according to
\begin{eqnarray*}
\;\oooo{1}{0}{1}{0}\;\otimes\;\oooo{1}{0}{0}{1}\;& =&\;\oooo{2}{0}{1}{1}\;\oplus\;\oooo{0}{1}{1}{1}\;\oplus\;\oooo{2}{1}{0}{0}\;\\
&&\oplus\;\oooo{1}{0}{0}{2}\;\oplus\;\oooo{0}{2}{0}{0}\;\\
&&\oplus 2\times\;\oooo{1}{0}{1}{0}\;
\;\oplus\;\oooo{0}{0}{0}{1},\;
\end{eqnarray*}
which compose as
$$\;\oooo{2}{0}{1}{1}\;=\;\xooo{2}{0}{1}{1}\;+
\begin{array}{c}
\;\xooo{1}{0}{0}{2}\;\\
\oplus\\
\;\xooo{0}{1}{1}{1}\;\\
\oplus\\
\;\xooo{1}{0}{1}{0}\;
\end{array}
+
\begin{array}{c}
\;\xooo{-1}{1}{0}{2}\;\\
\oplus\\
\;\xooo{-2}{2}{1}{1}\;\\
\oplus\\
\;\xooo{0}{0}{0}{1}\;\\
\oplus\\
\;\xooo{-1}{1}{1}{0}\;
\end{array}
+
\begin{array}{c}
\;\xooo{-2}{1}{0}{1}\;\\
\oplus\\
\;\xooo{-3}{2}{1}{0}\;\\
\oplus\\
\;\xooo{-3}{2}{0}{2}\;
\end{array}
+
\;\xooo{-4}{2}{0}{1}\;,$$

$$\;\oooo{0}{1}{1}{1}\;=\;\xooo{0}{1}{1}{1}\;+
\begin{array}{c}
\;\xooo{-1}{0}{2}{1}\;\\
\oplus\\
\;\xooo{-1}{1}{0}{2}\;\\
\oplus\\
\;\xooo{-1}{1}{1}{0}\;
\end{array}
+
\begin{array}{c}
\;\xooo{-2}{0}{1}{2}\;\\
\oplus\\
\;\xooo{-2}{1}{0}{1}\;\\
\oplus\\
\;\xooo{-2}{0}{2}{0}\;
\end{array}
+
\;\xooo{-3}{0}{1}{1}\;,$$

$$\;\oooo{2}{1}{0}{0}\;=\;\xooo{2}{1}{0}{0}\;+
\begin{array}{c}
\;\xooo{1}{0}{1}{0}\;\\
\oplus\\
\;\xooo{0}{2}{0}{0}\;
\end{array}
+
\begin{array}{c}
\;\xooo{-1}{1}{1}{0}\;\\
\oplus\\
\;\xooo{-2}{3}{0}{0}\;
\end{array}
+
\;\xooo{-3}{2}{1}{0}\;,$$

$$\;\oooo{1}{0}{0}{2}\;=\;\xooo{1}{0}{0}{2}\;+
\begin{array}{c}
\;\xooo{0}{0}{0}{1}\;\\
\oplus\\
\;\xooo{-1}{1}{0}{2}\;
\end{array}
+
\begin{array}{c}
\;\xooo{-1}{0}{0}{0}\;\\
\oplus\\
\;\xooo{-2}{1}{0}{1}\;
\end{array}
+
\;\xooo{0}{1}{0}{0}\;,$$

$$\;\oooo{1}{0}{1}{0}\;=\;\xooo{1}{0}{1}{0}\;+
\begin{array}{c}
\;\xooo{0}{0}{0}{1}\;\\
\oplus\\
\;\xooo{-1}{1}{1}{0}\;
\end{array}
+
\;\xooo{-2}{1}{0}{1}\;$$
and
$$\;\oooo{0}{2}{0}{0}\;=\;\xooo{0}{2}{0}{0}\;+
\;\xooo{-1}{1}{1}{0}\;
+
\;\xooo{-2}{0}{2}{0}\;,
$$

$$\;\oooo{0}{0}{0}{1}\;=\;\xooo{0}{0}{0}{1}\;+
\;\xooo{-1}{0}{0}{0}\;.
$$
There are 5 first order bilinear invariant differential pairings according to the projections onto (including the weights $k=1+v$ for vector fields of projective weight $v$ and $m=w-3$ for 2-forms of projective weight $w$, i.e.~we have to tensor by $\mathcal{O}(k-M)\otimes\mathcal{O}(m-M)=\mathcal{O}(v+w-4)$):
$$\;\xooo{v+w-4}{2}{0}{0}\;,\;\xooo{v+w-3}{0}{0}{2}\;,\;\xooo{v+w-4}{1}{1}{1}\;$$
and the two projections onto 
$$\;\xooo{v+w-3}{0}{1}{0}\;=\Omega^{2}(v+w)\;,$$
corresponding to
$$\;\ooo{0}{0}{1}\;\otimes\g_{1}\otimes \;\ooo{0}{1}{0}\;=2\times \;\ooo{0}{1}{0}\;\oplus\;\ooo{1}{1}{1}\;\oplus\;\ooo{2}{0}{0}\;\oplus\;\ooo{0}{0}{2}\;.$$
The concrete formulae for the two projections onto $\Omega^{2}(v+w)$ were given at the end of section~4.

\subsection{Weighted functions of excluded geometric weight}
Returning to Example (3) in 7.0.5, let us assume that the central character of $M_{\p}(\mathbb{V}_{0}(w-M))$ equals the central character of $M_{\p}(\mathbb{V}_{l}(w-M))$,
i.e.~$0\leq w=l-1\leq M-1$. This corresponds to an $l$-th order invariant differential operator 
$$D:\xooo{w}{0}{0}{0}\;...\;\oo{0}{0}\;\rightarrow \xoo{w-2l}{l}{0}\;...\oo{0}{0}\;...\;\oo{0}{0}.$$
Hence one can invariantly write $D(f)$, for $f\in\mathcal{O}(w)$.
Now we look at the $\p$-module
$$\tilde{\mathbb{V}}_{M,l}(\C)(w-M)=\mathbb{V}_{l}(w-M)+\mathbb{V}_{l+1}(w-M)+...+\mathbb{V}_{M}(w-M).$$
The central character of  $M_{\p}(\mathbb{V}_{l}(k-M))$ is different from the central character of all the other generalized Verma modules, because each $M_{\p}(\mathbb{V}_{s}(w-M))$ has the same central character as $M_{\p}(\mathbb{V}_{0}(w-M))$ if and only if~$w=s-1$. Therefore we can define an invariant differential mapping
$$\mathcal{O}(w)\stackrel{D}{\rightarrow}V_{l}(w-M)\rightarrow \tilde{V}_{M,l}(\C)(w-M)\hookrightarrow V_{M}(\C)(w-M).$$
The invariant pairings that we obtain via this construction do not involve derivatives of $f$ of order smaller than~$l$. This is confirmed by the formulae obtained earlier. 
\par
These considerations yield:
\begin{corollary}
If $M\geq\mathrm{max}_{i}\{a_{i}\}$ and $k$ does not equal one of the excluded weights up to order $M$ for $V=\;\xoo{k}{a_{1}}{a_{2}}\;...\;\oo{a_{n-2}}{a_{n-1}}\;$, then there is 
a one parameter family of invariant bilinear differential pairings of order $M$ between sections of $V$ and arbitrarily weighted functions onto every bundle that is induced by
an irreducible component of $\odot^{M}\g_{1}\otimes\;\oo{a_{1}}{a_{2}}\;...\;\oo{a_{n-2}}{a_{n-1}}\;$.
\end{corollary}

\subsubsection{Remark}
Using Pierie's formula, it is clear that the tensor product $\odot^{M}\g_{1}\otimes\;\oo{a_{1}}{a_{2}}\;...\;\oo{a_{n-2}}{a_{n-1}}\;$ does not have multiplicities.

\subsubsection{Example}
Let us analyze the example given in section 5, where we considered second order pairings $T\mathbb{CP}_{n}(v)\times\mathcal{O}(w)\rightarrow\Omega^{1}(v+w)$. For this purpose, we decompose
$$\odot^{2}\g_{1}\otimes\oo{0}{0}\;...\;\oo{0}{1}\;=\;\oo{2}{0}\;...\;\oo{0}{1}\oplus\;\oo{1}{0}\;...\;\oo{0}{0}.$$
Therefore if $v\not=-1,-(n+1)$ (for the other projection we also need to exclude $v=0$), then there should be a second order invariant differential pairing. This is true and the formula was given in section~5. Moreover one can clearly see which terms vanish in case the weight $w$ is excluded.

\section{Appendix}

\begin{definition}
\rm Let $\mathcal{M}$ be a (smooth) complex manifold and $V,W$ holomorphic vector bundles over $\mathcal{M}$ as in~3.6.  For every holomorphic vector bundle $U$ over $\mathcal{M}$ and every integer $k\in\N$ there exists the associated jet bundle $J^{k}U$ and for every $0\leq l\leq k$ there is a projection $\pi^{k}_{l}:J^{k}U\rightarrow J^{l}U$ (see~~\cite{kms}, p.~117,
definition~12.2).
If $\Lambda^{1}$ denotes the cotangent bundle on $\mathcal{M}$, then the projections can be put into an exact sequence
$$0\rightarrow \odot^{k} \Lambda^{1}\otimes U\rightarrow J^{k}U\rightarrow J^{k-1}U\rightarrow 0$$
as described in~~\cite{s},~p.~182. This exact sequence induces a filtration
$$J^{k}U=\sum_{l=0}^{k}\odot^{l}\Lambda^{1}\otimes U=U+\Lambda^{1}\otimes U+\odot^{2}\Lambda^{1}\otimes U+...+\odot^{k}\Lambda^{1}\otimes U$$
on the jet bundle.
\par
The mapping
$$\varphi_{M}=\oplus_{k+l=M}\pi_{k}^{M}\otimes\pi^{M}_{l}:J^{M}V\otimes J^{M}W\rightarrow\bigoplus_{k+l=M}J^{k}V\otimes J^{l}W$$
defines a canonical subbundle $B=\mathrm{ker}\; \varphi_{M}$ in $J^{M}V\otimes J^{M}W$, so that
$$J^{M}(V,W)=(J^{M}V\otimes J^{M}W)/\mathrm{ker}\; \varphi_{M}.$$
\end{definition}

\subsubsection{Remark}
It is easy to see that the vector bundle $J^{M}(V,W)$ has a filtration
$$J^{M}(V,W)=\sum_{k=0}^{M}\bigoplus_{l=0}^{k}\odot^{l}\Lambda^{1}\otimes V\otimes\odot^{k-l}\Lambda^{1}\otimes W,$$
which is equivalent to a series of exact sequences
$$0\rightarrow \bigoplus_{l=0}^{k}\odot^{l}\Lambda^{1}\otimes V\otimes\odot^{k-l}\Lambda^{1}\otimes W\stackrel{\iota}{\rightarrow} J^{k}(V,W)\rightarrow J^{k-1}(V,W)\rightarrow 0,$$
for $0\leq k\leq M$. The exact sequence
$$0\rightarrow \bigoplus_{l=0}^{M}\odot^{l}\Lambda^{1}\otimes V\otimes\odot^{M-l}\Lambda^{1}\otimes W\stackrel{\iota}{\rightarrow} J^{M}(V,W)\rightarrow J^{M-1}(V,W)\rightarrow 0$$
gives rise to a {\bf symbol} $\sigma=\phi\circ\iota$ for every homomorphism $\phi:J^{M}(V,W)\rightarrow E$, i.e. for every $M$-th order bilinear differential pairing.
\par
In the homogeneous case $J^{M}(V,W)$ is a homogeneous bundle with a $\p$-module structure on the fibre $J^{M}(\mathbb{V},\mathbb{W})$ that is induced by the $\p$-module structures of $J^{M}\mathbb{V}$ and~$J^{M}\mathbb{W}$.

\subsubsection{Remark}
It is possible to approach the theory of invariant differential pairings on homogeneous spaces from a completely algebraic point of view by considering $\mathfrak{U}(\g)$-modules that are dual to $J^{\infty}(\mathbb{V},\mathbb{W})$. These modules can be constructed by appropriately generalizing the {\bf bi-Verma} modules defined in ~\cite{d}, p.~6. Invariant differential pairings then correspond to so called {\bf singular vectors} in the sense of ~\cite{d} and it is possible to write down explicit formulas for infinitely many of them (given $V$ and $W$) if certain weights, that correspond to singular vectors in $M_{\p}(\mathbb{V})$ or $M_{\p}(\mathbb{W})$ and therefore to invariant differential operators, are excluded. In contrast to bi-Verma modules which are suitable for finding symmetric pairings between identical (line) bundles ($V=W$) and hence non-linear operators, this theory is valid for all pairings between arbitrary vector bundles $V$, $W$ and for any parabolic subalgebra $\p\subset\g$.
\par 
For the results proved and calculations performed in this article,  
however, this algebraic viewpoint is (at the moment) simply an alternative and  provides no particular advantage.

\par
\vspace{1cm}

I would like to thank my supervisor Prof.~Michael Eastwood for suggesting the problem and for his continuing help and support.

\end{document}